\documentclass[12 pt]{article}
\usepackage{amssymb,amsmath,amsfonts,amsthm}
\usepackage{graphicx}
\newtheorem {Theorem}{Theorem}[section]

\newtheorem {Proposition}[Theorem]{Proposition}

\theoremstyle{definition}
\newtheorem{Definition}{Definition}[section]

\newtheorem{Remark}{Remark}[section]

\newcommand\beq{\begin{equation}}
\newcommand\eeq{\end{equation}}

\begin{document}
\title{ The Mann-Whitney U-statistic for $\alpha$-dependent sequences}

\author{ J\'er\^ome Dedecker$^a$ and Guillaume Sauli\`ere$^b$}

\maketitle

{
\small{$^a$ Universit\'e Paris Descartes, Sorbonne Paris Cit\'e,  Laboratoire MAP5 (UMR 8145).\\
Email: jerome.dedecker@parisdescartes.fr 

\smallskip

$^b$ Universit\'e Paris Sud, Institut de Recherche bioM\'edicale et d'Epid\'emiologie du Sport (IRMES), Institut National du Sport de l'Expertise et de la Performance (INSEP).\\
Email:  Guillaume.SAULIERE@insep.fr }
} 

{\abstract{We give the asymptotic behavior of the Mann-Whitney U-statistic for two independent
stationary sequences.  The result applies to a large class of short-range dependent sequences,
including many non-mixing processes in the sense of Rosenblatt \cite{R}. We also give some 
partial results in the long-range dependent case, 
and we  investigate other related questions. Based on the theoretical results, 
we propose some simple corrections of the usual tests for stochastic domination; next we simulate different (non-mixing) stationary processes to see that the corrected tests perform well.}}

\section{Introduction and main result}
Let $(X_i)_{i \in {\mathbb Z}}$ and $(Y_i)_{i \in {\mathbb Z}}$
be two independent stationary sequences of real-valued random variables. Let $m=m_n$ be 
a sequence of positive integers such that $m_n \rightarrow \infty$ as $n \rightarrow \infty$. The Mann-Whitney U-statistic
may be defined as
$$
U_n= U_{n,m}= \frac{1}{n m } 
\sum_{i=1}^n \sum_{j=1}^m  {\bf 1}_{X_i < Y_j}  \, .
$$
We shall also use the notations $\pi= {\mathbb P}(X_1 < Y_1)$, 
$H_Y(t)= {\mathbb P}(Y_1>t)$ and 
$G_X(t)={\mathbb P}(X_1<t)$.  

It is well known  that, 
if $(X_i)_{i \in {\mathbb Z}}$ and $(Y_i)_{i \in {\mathbb Z}}$ are 
two independent sequences of independent and identically (iid)
random variables, then $U_n$ may be used to test 
$\pi=1/2$ against $\pi \neq 1/2$ (or $\pi >1/2$, $\pi < 1/2$).
When $\pi >1/2$ we shall say that $Y$ stochastically dominates 
$X$ in a weak sense. This property of weak domination  is true for instance
 if $Y$ stochastically dominates $X$, that is if $H_Y(t) \geq H_X(t)$ for any real $t$, with a strict inequality for some  $t_0$.
 But  it also holds in many other situations, for instance if
 $X$ and $Y$ are two Gaussian random variables with 
 ${\mathbb E}(X)< {\mathbb E}(Y)$, whatever the variances
 of $X$ and $Y$.
 
 In this paper, we study the asymptotic behavior of $U_n$ under some conditions on the $\alpha$-dependence coefficients  of 
 the sequences $(X_i)_{i \in {\mathbb Z}}$ and $(Y_i)_{i \in {\mathbb Z}}$. Let us recall the definition of these coefficients.
 
 \begin{Definition}
 Let $(X_i)_{i \in {\mathbb Z}}$ be a strictly stationary sequence, 
 and let ${\mathcal F}_0=\sigma (X_k, k \leq 0)$. For any non-negative integer $n$, let 
 $$
 \alpha_{1,\bf{X}}(n)=
 \sup_{x \in {\mathbb R}}\left \|{\mathbb E}\left ({\bf 1}_{X_n \leq x}|{\mathcal F}_0\right ) - F(x)
\right \|_1 \, ,
$$ 
where $F$ is the distribution function of $X_0$. 
 \end{Definition}
 
 We shall also always require that the sequences are $2$-ergodic in the following sense.
 
 \begin{Definition}
 A strictly stationary sequence $(X_i)_{i \in {\mathbb Z}}$ is 
 $2$-ergodic if, for any bounded mesurable function $f$ from 
 ${\mathbb R}^2$ to ${\mathbb R}$,  and any non-negative integer $k$
 \begin{equation}\label{2erg}
 \lim_{n \rightarrow \infty} \frac 1 n \sum_{i=1}^n f(X_i, X_{i+k}) =
 {\mathbb E}(f(X_0, X_k)) \quad \text{almost surely.}
 \end{equation}
 \end{Definition}
 Note that the almost sure  limit  in \eqref{2erg} always exists, by the ergodic
 theorem. The property of 2-ergodicity means only that this 
 limit is constant. This property is weaker than usual ergodicity, 
 which would implies the same property for  any bounded function $f$
 from ${\mathbb R}^\ell$ to ${\mathbb R}$, $\ell \in {\mathbb N}-\{0\}$. 
 
 Our main result is the following theorem.

\begin{Theorem}\label{MannTh}
Let $(X_i)_{i \in {\mathbb Z}}$ and $(Y_i)_{i \in {\mathbb Z}}$
be two independent stationary 2-ergodic sequences
of real-valued random variables. Assume that 
\begin{equation}\label{alphacond}
\sum_{k=1}^\infty \alpha_{1,\bf{X}}(k) < \infty  \quad 
\text{and} \quad 
\sum_{k=1}^\infty \alpha_{1,\bf{Y}}(k) < \infty \, .
\end{equation} 
Let $(m_n)$ be a sequence of positive integers such that
\begin{equation}\label{cond1}
\quad \lim_{n \rightarrow \infty} \frac{n}{m_n} = c
\, ,
\quad \text{for some $c \in [0, \infty)$.}
\end{equation}
Then the random variables 
$
  \sqrt{n} (U_n- \pi) 
$
converge in distribution to ${\mathcal N}(0,V)$, with
\begin{multline}\label{Var}
V=\mathrm{Var}(H_Y(X_0)) + 2 \sum_{k=1}^\infty 
\mathrm{Cov}(H_Y(X_0), H_Y(X_k)) \\ + 
c \left ( \mathrm{Var}(G_X(Y_0)) + 2 \sum_{k=1}^\infty 
\mathrm{Cov}(G_X(Y_0), G_X(Y_k))\right ) \, .
\end{multline}
\end{Theorem}

 Note that the result of  Theorem \ref{MannTh} has been obtained by Serfling \cite{S}, 
provided that the more restrictive condition
\begin{equation}
 \sum_{k \geq 0} \alpha(k) < \infty \quad \text{and} \quad \alpha(k)= O\left( \frac{1}{k \log k} \right) 
\end{equation}
is satisfied by the two sequences $(X_i)_{i \in {\mathbb Z}}$ and $(Y_i)_{i \in {\mathbb Z}}$, 
where the $\alpha(k)$'s are  the usual $\alpha$-mixing coefficients of  Rosenblatt \cite{R}.

\smallskip

The notion of $\alpha$-dependence that we consider in the present paper is  much weaker 
than  $\alpha$-mixing in the sense of Rosenblatt \cite{R}. 
Some properties of these $\alpha$-dependence coefficients are stated in the book by Rio 
\cite{Ri}, and  many examples are  given in the two papers 
\cite{DP1}, \cite{DP2}. Let us recall one of these 
examples. Let 
\begin{equation}\label{proclin}
 X_i= \sum_{k \geq 0} a_k \varepsilon_{i-k}
\end{equation}
where  $(a_k)_{k \geq 0} \in \ell^2$, and the sequence $(\varepsilon_i)_{i \in {\mathbb Z}}$
is a sequence of independent and identically distributed 
real-valued random variables, with mean zero and finite variance.
If the cumulative distribution function $F$ of $X_0$ is $\gamma$-H\"older 
for some $\gamma \in (0, 1]$, then 
$$
\alpha_{1,\bf{X}}(n) \leq  C\left( \sum_{k \geq n} a_k^2 \right)^{\frac{\gamma}{\gamma +2}}
$$
for some positive constant $C$ (this follows from the computations in 
Section 6.1 of \cite{DP2}). In particular, if $a_k$ is geometrically decreasing, then so is $\alpha_{1,\bf{X}}(n)$ (whatever the index $\gamma$). 
 
 Note that, without extra assumptions on the distribution of $\varepsilon_0$, such linear processes have no reasons to be mixing in the sense of 
Rosenblatt. For instance, it is well known that the linear process
$$
X_i= \sum_{k \geq 0} \frac{\varepsilon_{i-k} }{2^{k+1}}, \quad
\text{where ${\mathbb P}(\varepsilon_1=-1/2)={\mathbb P}(\varepsilon_1=1/2)=1/2$,}
$$
 is not $\alpha$-mixing (see
for instance \cite{Br}). By contrast, for this particular 
example, $\alpha_{1,\bf{X}}(n)$ is geometrically decreasing.

To conclude, let us also mention the paper \cite{DGM} where it is proved that the coefficients $\alpha_{1,\bf{X}}(n)$ can be computed for a large class of dynamical systems on the unit interval
(see Subsection \ref{Sec:IM} of the present paper for more details). 

\begin{Remark}
Other types of dependence are considered in the papers by 
Dewan and Prakasa Rao \cite{DPR} and Dehling and Fried \cite{DF}. The paper
\cite{DPR}  deals with the U-statistics 
$U_n$ for two independent samples of associated random variables. The second paper \cite{DF}
deals with general U-statistics of functions of $\beta$-mixing sequences, for either 
two independent samples, or for two adjacent samples from the same time series
(see Subsection \ref{adj} for an extension of our results to this more difficult context). 
The class of unctions of $\beta$-mixing sequences contains also many  examples,
with a large intersection with the class of $\alpha$-dependent sequences. The 
advantage of our approach is that it leads to conditions that are optimal
in a precise  sense (see Section \ref{Sec:LRD}).  On another hand, it seems more difficult to deal with
general U-statistics in our context (some kind of regularity is needed, for instance we 
could certainly deal with U-statistics  based on functions which have bounded variation in each coordinates).

\end{Remark}

\section{Testing the stochastic domination from two independent time series}\label{sec2}

\setcounter{equation}{0}

The statistic $\sqrt n (U_n-1/2)$ cannot be used to test $H_0: \pi=1/2$
versus one of the possible alternatives, because under $H_0$
the asymptotic distribution of $\sqrt n (U_n-1/2)$ depends of  the unknow 
quantity $V$ defined in \eqref{Var}. In the next  Proposition, we propose an 
estimator of  $V$ under some slightly 
stronger conditions on the dependence coefficients. 

\begin{Definition}\label{defalpha}
 Let $(X_i)_{i \in {\mathbb Z}}$ be a strictly stationary sequence
 of real-valued random variables, 
 and let ${\mathcal F}_0=\sigma (X_k, k \leq 0)$. 
 Let $f_x(z)= {\bf 1}_{z \leq x}-F(x)$, where $F$ is the cumulative distribution function of  $X_0$. For any non-negative integer $n$, let
\begin{equation*}
 \alpha_{2,\bf{X}}(n)= \max \left \{ \alpha_{1,\bf{X}}(n), 
 \sup_{x, y \in {\mathbb R}, j\geq i \geq n}\left \|{\mathbb E}\left (f_x(X_i)f_y(X_j)|{\mathcal F}_0\right ) - {\mathbb E}\left (  f_x(X_i)f_y(X_j)\right ) 
\right \|_1  \right \}\, .
\end{equation*}
 \end{Definition}

\begin{Proposition}\label{varest}
Let $(X_i)_{i \in {\mathbb Z}}$ and $(Y_i)_{i \in {\mathbb Z}}$
be two independent stationary sequences, and let 
$(m_n)$ be a sequence of positive inetegers satisfying
\eqref{cond1}.  Assume that 
\begin{equation}\label{alphacond2}
\sum_{k=1}^\infty \alpha_{2,\bf{X}}(k) < \infty  \quad 
\text{and} \quad 
\sum_{k=1}^\infty \alpha_{2,\bf{Y}}(k) < \infty \, .
\end{equation} 
Let 
$$
H_m(t)= \frac 1 m \sum_{j=1}^m {\bf 1}_{Y_j>t} \quad \text{and}
\quad 
G_n(t)= \frac 1n \sum_{i=1}^n {\bf 1}_{X_i<t} \, , 
$$
and 
$$
 \hat \gamma_X (k)= \frac 1 n  \sum_{i=1}^{n-k}
 (H_m(X_i)-\bar{H}_m )(H_m(X_{i+k}) - \bar{H}_m) \, ,\quad 
  \hat \gamma_Y (\ell)= \frac 1 m  \sum_{j=1}^{m-\ell}
 (G_n(Y_j)-\bar{G}_n )(G_n(Y_{j+\ell}) - \bar{G}_n) \, ,
$$
where $n \bar{H}_m=  \sum_{i=1}^{n}
 H_m(X_i)$ and $m \bar{G}_n=  \sum_{j=1}^{m}
 G_n(Y_j)$. 
Let also $(a_n)$ and $(b_n)$ be two sequences of positive
 integers tending to infinity as $n$ tends to infinity, such that $a_n=o(\sqrt n)$ and $b_n=o(\sqrt {m_n})$. 
 Then 
 $$
 V_n=  \hat \gamma_X (0) +2 \sum_{k=1}^{a_n} \hat \gamma_X (k)  
 +\frac{n}{m_n} \left ( 
 \hat \gamma_Y (0) +2 \sum_{\ell=1}^{b_n} \hat \gamma_Y (\ell)   \right )
 $$
 converges in ${\mathbb L}^2$ to the quantity $V$ defined
 in \eqref{Var}.
\end{Proposition}

Combining Theorem \ref{MannTh} and Proposition \ref{varest}, 
we obtain that, under $H_0: \pi=1/2$, if $V>0$, the random variables
\begin{equation}\label{Tn}
T_n= \frac{\sqrt{n}(U_n - 1/2)}{\sqrt{\max \{V_n, 0\}}}
\end{equation}
converge in distribution to ${\mathcal N}(0,1)$. 

\begin{Remark}
The condition \eqref{alphacond2} is not restrictive: in all the 
natural examples we know, the coefficients $\alpha_{2, \bf X}(k)$
behave as $\alpha_{1, \bf X}(k)$ (this is the case for instance for 
the linear process \eqref{proclin}). Note that, if we assume \eqref{alphacond2} instead of \eqref{alphacond}, the ergodicity is no longer required in Theorem 
\ref{MannTh}. Finally, 
following the proof of Theorem 4.1 page 89 in \cite{De}, one can
also build a consistent  estimator of $V$ under the condition 
\begin{equation}\label{alphacond3}
\sum_{k=1}^\infty \sqrt{\frac{\alpha_{1,\bf{X}}(k)}{k} }< \infty  \quad 
\text{and} \quad 
\sum_{k=1}^\infty \sqrt{\frac{\alpha_{1,\bf{Y}}(k)}{k} }< \infty 
 \, ,
\end{equation} 
which is not directly comparable to \eqref{alphacond2} 
(for instance, the first part of  \eqref{alphacond2}  is satisfied  if 
$\alpha_{2,\bf{X}}(k)= O(k^{-1} (\log k)^{-a})$ for some $a>1$, while 
the first part of \eqref{alphacond3} requires 
$\alpha_{1,\bf{X}}(k)= O(k^{-1} (\log k)^{-a})$ for some $a>2$). 
\end{Remark}

\begin{Remark}
The choice of the two sequences $(a_n)$ and $(b_n)$ is a delicate
matter. If the coefficients $\alpha_{2,\bf{X}}(k)$ decrease very 
quickly, then $a_n$ should increase very slowly (it suffices to take
$a_n \equiv 0$ in the iid setting). On the contrary, if 
$\alpha_{2,\bf{X}}(k)= O(k^{-1} (\log k)^{-a})$ for some 
$a>1$, then the terms in the covariance series have no 
reason to be small, and one should take $a_n$ close 
to $\sqrt n$ to estimate many of these covariance terms.
A data-driven criterion for choosing $a_n$ and $b_n$ is 
an interesting (but probably difficult) question, which is far beyond the 
scope of the present paper. 

However, from a practical point of view, there is an easy way 
to proceed: one can plot the estimated covariances $\hat \gamma_X(k)$'s and choose
$a_n$ (not too large) in such a way that 
$$
\hat \gamma_X (0) +2 \sum_{k=1}^{a_n} \hat \gamma_X (k)
$$
should represent an important part of the covariance series
$$
\mathrm{Var}(H_Y(X_0)) + 2 \sum_{k=1}^\infty 
\mathrm{Cov}(H_Y(X_0), H_Y(X_k)) \, .
$$
Similarly, one may choose $b_n$ from the graph of 
the $\hat \gamma_Y(\ell)$'s. As we shall see in the simulations
(Section \ref{Sec:Simu}),
if the decays of the true covariances $ \gamma_X(k)$ and 
$\gamma_Y(\ell)$ are not too slow, this provides an easy and reasonable choice for
 $a_n$ and $b_n$.

\end{Remark}

\section{Long range dependence}\label{Sec:LRD}
\setcounter{equation}{0}
In this section, we shall see that the condition \eqref{alphacond}
cannot be weakened for the validity of Theorem \ref{MannTh}.
More precisely, we shall give some examples where 
\eqref{alphacond} fails, and the asymptotic 
distribution of $a_n (U_n-\pi)$ (if any) can be non Gaussian. 

Let us begin with the boundary case, where 
\begin{equation}\label{boundarycond}
\alpha_{2,\bf{X}}(k)= O((k+1)^{-1}) \quad  \text{and} \quad 
\alpha_{2,\bf{Y}}(k) = O((k+1)^{-1}) \, .
\end{equation}
In that case, one can prove 
 the following result.

\begin{Proposition}\label{BC}
Let $(X_i)_{i \in {\mathbb Z}}$ and $(Y_i)_{i \in {\mathbb Z}}$
be two independent stationary sequences. Assume that 
\eqref{boundarycond} is satisfied, and let $(m_n)$ be a sequence 
a positive integers satisfying \eqref{cond1}. Then:
\begin{enumerate}
\item For any $r \in (2,4)$, there exists a positive constant $C_r$ such that, for any $x>0$,
$$
  \limsup_{n \rightarrow \infty}
  {\mathbb P} \left ( \sqrt{n /\log n} (U_n-\pi)>x \right)
  \leq \frac{C_r}{x^r} \, .
$$
\item There are some examples of stationary Markov chains
$(X_i)_{i \in {\mathbb Z}}$ and $(Y_i)_{i \in {\mathbb Z}}$ for which the sequence $\sqrt{n/\log n} (U_n-\pi)$ converges in distribution 
to a non-degenerate normal distribution.
\end{enumerate}
\end{Proposition}

We consider now the case where 
\begin{equation}\label{LRcond}
 \quad \sum_{k=1}^\infty \alpha_{1,\bf{X}}(k) < \infty  \quad 
\text{and}   \quad 
\alpha_{1,\bf{Y}}(k) = O((k+1)^{-(p-1)}) \, , \ \text{for some $p \in (1,2)$}.
\end{equation}
This means that $(X_i)_{i \in {\mathbb Z}}$ is short-range dependent, while  $(Y_i)_{i \in {\mathbb Z}}$ is possibly 
long-range dependent. 
In that case, we shall prove the following result.

\begin{Proposition}\label{LR}
Let $(X_i)_{i \in {\mathbb Z}}$ and $(Y_i)_{i \in {\mathbb Z}}$
be two independent stationary sequences. Assume that 
\eqref{LRcond} is satisfied for some $p \in (1,2)$, 
and that $(X_i)_{i \in {\mathbb Z}}$ is 2-ergodic. Let $(m_n)$ be a sequence 
a positive integers such that 
\begin{equation}\label{cond2}
\quad \lim_{n \rightarrow \infty} \frac{n}{m_n^{2(p-1)/p}} = c
\, ,
\quad \text{for some $c \in [0, \infty)$.}
\end{equation} Then:
\begin{enumerate}
\item There exists a positive constant $C_p$ such that, for any $x>0$,
$$
  \limsup_{n \rightarrow \infty}
  {\mathbb P} \left ( \sqrt n (U_n-\pi)>x \right)
  \leq \frac{C_p}{x^p} \, .
$$
\item There are some examples of stationary Markov chains
$(X_i)_{i \in {\mathbb Z}}$ and $(Y_i)_{i \in {\mathbb Z}}$ for which the sequence $\sqrt{n} (U_n-\pi)$ converges in distribution 
to $Z+\sqrt c S$, where $Z$ is independent of $S$, $Z$ is a Gaussian random variable, and $S$ is a stable random variable of order $p$.
\end{enumerate}
\end{Proposition}

\begin{Remark}
Note that we do not deal with the case where the two samples are long-range dependent;
in that case the degenerate U-statistic from the Hoeffding decomposition (see \eqref{Hoef})  is no longer negligible, 
which gives rise to other possible limits. See the paper \cite{DRT} for a treatment of this 
difficult question 
in the case of functions of Gaussian processes. 
\end{Remark}

\section{Proofs}
\setcounter{equation}{0}
\subsection{Proof of Theorem \ref{MannTh}}

The main ingredient of the proof is the following proposition:

\begin{Proposition}\label{degen}
Let $(X_i)_{i \in {\mathbb Z}}$ and $(Y_i)_{i \in {\mathbb Z}}$
be two independent stationary sequences. Let also 
\begin{equation}\label{degenU}
f(x,y)={\bf 1}_{x<y}-H_Y(x)-G_X(y) +  \pi \, .
\end{equation}
Then
$$
{\mathbb E}\left (\left( \frac{1}{nm} \sum_{i=1}^n \sum_{j=1}^m
f(X_i,Y_j) \right )^2 \right ) \leq \frac{16}{mn}
\sum_{i=0}^{n-1} \sum_{j=0}^{m-1} \min \{ \alpha_{1, \bf X}(i), \alpha_{1, \bf Y}(j) \} \, .
$$
\end{Proposition}

Let us admit this proposition for a while, and see how it can be used to prove the theorem. Starting from \eqref{degenU}, we easily see that
\begin{equation}\label{Hoef}
\sqrt n (U_n-\pi)= \frac{\sqrt n}{nm} \sum_{i=1}^n \sum_{j=1}^m
f(X_i,Y_j) + \frac{1}{ \sqrt n} \sum_{i=1}^n (H_Y(X_i)- \pi) 
+ 
\frac {\sqrt n}{m} \sum_{j=1}^m (G_X(Y_j)- \pi) \, .
\end{equation}
This is the well known Hoeffding decomposition of U-statistics; 
the first term  on the right-hand side of \eqref{Hoef} is 
 called a \emph{degenerate} U-statistic, and is handled \emph{via} Proposition \ref{degen}. Indeed, it follows easily from this proposition that 
\begin{equation}\label{secondbound}
n\mathbb{E}\left ( \left( \frac{1}{nm} \sum_{i=1}^n \sum_{j=1}^m
f(X_i,Y_j) \right )^2 \right ) \leq \frac{16}{m} \sum_{i=0}^{n-1} (i+1) 
\alpha_{1, \bf X}(i) + 
\frac{16}{m} \sum_{j=0}^{m-1} (j+1) 
\alpha_{1, \bf Y}(j) \, .
\end{equation}
Since \eqref{alphacond} holds, then $\alpha_{1, \bf X}(n)=o(n^{-1})$ 
and  $\alpha_{1, \bf Y}(n)=o(m^{-1})$ (because these coefficients 
decrease) and by Cesaro's lemma
$$
\lim_{ n \rightarrow \infty}\frac 1 n \sum_{i=0}^n (i+1) \alpha_{1, \bf X}(i)=0 \quad  \text{and} \quad 
\lim_{ m \rightarrow \infty}\frac 1 m \sum_{j=0}^m (j+1) \alpha_{1, \bf Y}(j)=0 \, .
$$
From \eqref{secondbound} and Condition \eqref{cond1}, 
we infer that
\begin{equation}\label{one}
\frac{\sqrt n}{nm} \sum_{i=1}^n \sum_{j=1}^m
f(X_i,Y_j) \quad \text{converges to 0 in ${\mathbb L}^2$ as $n \rightarrow \infty$.} 
\end{equation}
 It remains to control the two last terms on the right-hand side of \eqref{Hoef}. We first note that these two terms are independent, so it is enough to prove  the convergence in distribution for both terms.
The functions $H_Y$ and $G_X$ being monotonic, the central
limit theorem for stationary and 2-ergodic $\alpha$-dependent
sequences  (see for instance \cite{DGM}) implies that
\begin{multline}\label{two}
\frac {1}{ \sqrt n} \sum_{i=1}^n (H_Y(X_i)- \pi)  \quad \text{converges in distribution to ${\mathcal N}(0, v_1)$, where}
\\ v_1= \mathrm{Var}(H_Y(X_0)) + 2 \sum_{k=1}^\infty 
\mathrm{Cov}(H_Y(X_0), H_Y(X_k))\, ,
\end{multline}
and 
\begin{multline}\label{three}
\frac {1}{ \sqrt m} \sum_{j=1}^m (G_X(Y_j)- \pi)  \quad \text{converges in distribution to ${\mathcal N}(0, v_2)$, where}
\\ v_2= \mathrm{Var}(G_X(Y_0)) + 2 \sum_{k=1}^\infty 
\mathrm{Cov}(G_X(Y_0), G_X(Y_k))\, .
\end{multline}
Combining \eqref{Hoef}, \eqref{one}, \eqref{two} and \eqref{three}, the proof of Theorem \ref{MannTh} is complete.

It remains to prove Proposition \ref{degen}. We start 
from the elementary inequality
\begin{equation} \label{elem}
{\mathbb E}\left (\left( \frac{1}{nm} \sum_{i=1}^n \sum_{j=1}^m
f(X_i,Y_j) \right )^2 \right )
\leq \frac{4}{n^2 m^2} \sum_{i=1}^n \sum_{j=1}^m
\sum_{k=i}^n \sum_{\ell =j}^m 
{\mathbb E}( f(X_i, Y_j) f(X_k, Y_\ell)) \, .
\end{equation}
To control the term ${\mathbb E}( f(X_i, Y_j) f(X_k, Y_\ell))$, 
we first
work conditionally on ${\bf Y}=(Y_i)_{i \in {\mathbb Z}}$. Note first that, since ${\bf Y}$ is independent of ${\bf X}$, ${\mathbb E}(f(X_i,Y_j)|{\bf Y})=0$. 
Since $|f(x,y)|\leq 2$, 
\begin{equation}\label{cov1}
|{\mathbb E}( f(X_i, Y_j) f(X_k, Y_l)| {\bf Y}={\bf y})|
\leq 2 \|{\mathbb E}(f(X_k, y_\ell) |{\mathcal F}_i)\|_1 \, ,
\end{equation}
where ${\mathcal F}_i= \sigma (X_k, k \leq i)$. Now, by definition of $f$, 
\begin{equation}\label{cov2}
\|{\mathbb E}(f(X_k, y_\ell) |{\mathcal F}_i)\|_1 \leq 
\|{\mathbb E}({\bf 1}_{X_k < y_\ell}|{\mathcal F}_i)-G_X(y_\ell)\|_1 +
\|{\mathbb E}(H_Y(X_k) |{\mathcal F}_i)- \pi \|_1 \, 
\leq 2 \alpha_{1, {\bf X}}(k-i) \, .
\end{equation}
From \eqref{cov1} and \eqref{cov2}, it follows that 
\begin{equation} \label{cov3}
|{\mathbb E}( f(X_i, Y_j) f(X_k, Y_\ell))| 
\leq 
\|{\mathbb E}( f(X_i, Y_j) f(X_k, Y_\ell)| {\bf Y})\|_1
\leq 
 4 \alpha_{1, {\bf X}}(k-i) \, .
\end{equation}
In the same way
\begin{equation}\label{cov4}
|{\mathbb E}( f(X_i, Y_j) f(X_k, Y_\ell))| 
\leq 
 4 \alpha_{1, {\bf Y}}(\ell-j) \, .
\end{equation}
Proposition \ref{degen} follows from \eqref{elem}, \eqref{cov3} and
\eqref{cov4}.

\subsection{Proof of Proposition \ref{varest}}
We proceed in two steps.

\smallskip

\emph{Step 1.} Let
$$
  \gamma^*_X (k)= \frac 1 n  \sum_{i=1}^{n-k}
 (H_Y(X_i)-\pi )(H_Y(X_{i+k}) - \pi) \, ,\quad 
  \gamma^*_Y (\ell)= \frac 1 m  \sum_{j=1}^{m-\ell}
 (G_X(Y_j)-\pi )(G_X(Y_{i+k}) - \pi) \, ,
$$
and note that 
\begin{equation}\label{esp}
{\mathbb E}(\gamma^*_X (k))=\frac{n-k}{n}\mathrm{Cov}(H_Y(X_0),H_Y(X_k)) \quad \text{and} \quad 
{\mathbb E}(\gamma^*_Y (\ell))=\frac{n-k}{n}\mathrm{Cov}(G_X(Y_0),G_X(Y_\ell))\, .
\end{equation}

In this first step,  we shall   prove that 
$$
 V^*_n= \gamma^*_X (0) +2\sum_{k=1}^{a_n}  \gamma^*_X (k)  
 +\frac{n}{m_n} \left ( \gamma^*_Y (0) +2 \sum_{\ell=1}^{b_n}  \gamma^*_Y (\ell)   \right ) \, ,
 $$ converges in ${\mathbb L}^2$ to $V$. Clearly, by \eqref{esp},  it suffices 
 to show that $\mathrm{Var}(V^*_n)$ converges to 0 as
 $n \rightarrow \infty$. We only deal with the second term
 in the expression of $V_n^*$, the other ones may be handled
 in the same way. Letting $Z_i=H_Y(X_i)-\pi$, 
 and $\gamma_k= \mathrm{Cov}(H_Y(X_0),H_Y(X_k))$ we have
 \begin{equation}\label{mainvar}
\mathrm{Var} \left ( \sum_{k=1}^{a_n}  \gamma^*_X (k)  \right) 
= \frac{1}{n^2}
 \sum_{k=1}^{a_n} \sum_{\ell=1}^{a_n} 
\sum_{i=1}^{n-k} \sum_{j=1}^{n-\ell}
 {\mathbb E}((Z_iZ_{i+k}-\gamma_k) (Z_jZ_{j+\ell}-\gamma_\ell) ) \, .
\end{equation}
 We must examine several cases.
 
 On $\Gamma_1=\{j \geq i+k\} $,
 $$
 |{\mathbb E}((Z_iZ_{i+k}-\gamma_k) (Z_jZ_{j+\ell}-\gamma_\ell) )| \leq  2 \alpha_{2, \bf X}(j-i-k).
 $$
  Consequently
 \begin{equation}\label{C1}
 \frac{1}{n^2}
 \sum_{k=1}^{a_n} \sum_{\ell=1}^{a_n} 
\sum_{i=1}^{n-k} \sum_{j=1}^{n-\ell}
|
 {\mathbb E}((Z_iZ_{i+k}-\gamma_k) 
 (Z_jZ_{j+\ell}-\gamma_\ell) )
|
 {\bf 1}_{\Gamma_1} \leq 
 \frac{2}{n^2} \sum_{k=1}^{a_n} \sum_{\ell=1}^{a_n} 
\sum_{i=1}^{n} \left ( \sum_{j \geq 0} \alpha_{2, \bf X}(j) \right)
 \leq \frac{C a_n^2}{n}
 \end{equation}
 for some positive constant $C$. 
 
 On $\Gamma_2=\{i \leq j < i+k \} \cap \{i+k \leq j +\ell \}$,
 $$
 |{\mathbb E}((Z_iZ_{i+k}-\gamma_k) (Z_jZ_{j+\ell}-\gamma_\ell) )|= |{\mathbb E}((Z_iZ_{i+k}-\gamma_k) Z_jZ_{j+\ell} )| \leq  2 \alpha_{1, \bf X}(j + \ell-i-k).
 $$
  Consequently
 \begin{equation}\label{C2}
 \frac{1}{n^2}
 \sum_{k=1}^{a_n} \sum_{\ell=1}^{a_n} 
\sum_{i=1}^{n-k} \sum_{j=1}^{n-\ell}
| {\mathbb E}((Z_iZ_{i+k}-\gamma_k)
 (Z_jZ_{j+\ell}-\gamma_\ell) ) |
 {\bf 1}_{\Gamma_2} \leq 
 \frac{2}{n^2} \sum_{k=1}^{a_n} \sum_{\ell=1}^{a_n} 
\sum_{i=1}^{n} \left ( \sum_{j \geq 0} \alpha_{1, \bf X}(j) \right)
 \leq \frac{C a_n^2}{n} \, .
 \end{equation}
 
 On $\Gamma_3=\{i \leq j  \} \cap \{i+k > j +\ell \}$,
 $$
 |{\mathbb E}((Z_iZ_{i+k}-\gamma_k) (Z_jZ_{j+\ell}-\gamma_\ell) )|= |{\mathbb E}((Z_jZ_{j+\ell}-\gamma_\ell) Z_iZ_{i+k} )| \leq  2 \alpha_{1, \bf X}(i+k-j - \ell).
 $$
 Consequently
 \begin{equation}\label{C3}
 \frac{1}{n^2}
 \sum_{k=1}^{a_n} \sum_{\ell=1}^{a_n} 
\sum_{i=1}^{n-k} \sum_{j=1}^{n-\ell}
 |{\mathbb E}((Z_iZ_{i+k}-\gamma_k) (Z_jZ_{j+\ell}-\gamma_\ell) )|
 {\bf 1}_{\Gamma_3} \leq 
 \frac{2}{n^2} \sum_{k=1}^{a_n} \sum_{\ell=1}^{a_n} 
\sum_{j=1}^{n} \left ( \sum_{i \geq 0} \alpha_{1, \bf X}(i) \right)
 \leq \frac{C a_n^2}{n} \, .
 \end{equation}
 
 From \eqref{C1}, \eqref{C2} and \eqref{C3}, setting
 $\Gamma=\{i\leq j\}= \Gamma_1 
 \cup \Gamma_2 \cup \Gamma_3$, one gets
 $$
 \frac{1}{n^2}
 \sum_{k=1}^{a_n} \sum_{\ell=1}^{a_n} 
\sum_{i=1}^{n-k} \sum_{j=1}^{n-\ell}
 | {\mathbb E}((Z_iZ_{i+k}-\gamma_k) 
 (Z_jZ_{j+\ell}-\gamma_\ell) )|
 {\bf 1}_{\Gamma} 
 \leq \frac{3C  a_n^2}{n} 
 $$
 Of course, interchanging $i$ and $j$, the same is true 
 on $\Gamma^c=\{i> j\}$, and finally
 \begin{equation}\label{C4}
 \frac{1}{n^2}
 \sum_{k=1}^{a_n} \sum_{\ell=1}^{a_n} 
\sum_{i=1}^{n-k} \sum_{j=1}^{n-\ell}
| {\mathbb E}((Z_iZ_{i+k}-\gamma_k) 
 (Z_jZ_{j+\ell}-\gamma_\ell) )|
 \leq \frac{6C  a_n^2}{n} 
 \end{equation}
 
 From \eqref{mainvar}, \eqref{C4} and the fact that $a_n=o(\sqrt n)$, we infer that 
 $\sum_{k=1}^{a_n}  \gamma^*_X (k)$ converges in ${\mathbb 
 L}^2$ to  $\sum_{k>0} \mathrm{Cov}(H_Y(X_0),H_Y(X_k))$. 
 Since the other terms in the definition of $V_n^*$ can be handled in the same way, it follows that $V_n^*$ converges 
in ${\mathbb L}^2$ to $V$. 

\smallskip

\emph{Step 2.} In this second step, in the expression of $\gamma_X^*(k)$, we replace 
$H_Y(X_i)$ by $H_m(X_i)$,  $H_Y(X_{i+k})$ by $H_m(X_{i+k})$, and $\pi$ by $\bar H_m$. Similarly, in the expression of $\gamma_Y^*(\ell)$, we replace
$G_X(Y_j)$ by $G_n(Y_j)$, $G_X(Y_{j+ \ell})$ by $G_n(Y_{j+ \ell})$, and $\pi$ by $\bar G_n$. If we can prove that these replacements in the expression of $V_n^*$ are negligible in ${\mathbb L}^2$, 
the conclusion of Proposition \ref{varest} will follow. 

Let 
$$
\gamma_{1,X} (k)= \frac 1 n  \sum_{i=1}^{n-k}
 (H_m(X_i)-\pi )(H_Y(X_{i+k}) - \pi)  \, .
$$
Then
\begin{equation}\label{rep1}
 \|\gamma_X^*(k)-\gamma_{1,X} (k)\|_2 \leq \| H_m(X_i)-H_Y(X_i)\|_2 = \left (\int \| H_m(x)-H_Y(x)\|_2^2 {\mathbb P}_X(dx) \right )^{1/2} \, ,
\end{equation}
the last equality being true because $(X_i)_{i \in {\mathbb Z}}$ 
and $(Y_i)_{i \in {\mathbb Z}}$ are independent. Now,
\begin{equation}\label{emp1}
\| H_m(x)-H_Y(x)\|_2^2 \leq \frac{1}{ m} \left ( 
\mathrm{Var}({\bf 1}_{Y_0>x}) +
2  \sum_{k=1}^{n-1} |\mathrm{Cov}({\bf 1}_{Y_0>x}, 
{\bf 1}_{Y_k>x})| \right ) \leq \frac{ 2 \sum_{k = 0}^{n-1}
\alpha_{1, \bf Y}(k)}{m} \, .
\end{equation}
From \eqref{rep1} and \eqref{emp1}, we infer that there exist a 
positive constant $C_1$ such that 
$$
\left \|\gamma^*_X (0)-\gamma_{1,X} (0) +2\sum_{k=1}^{a_n}  (\gamma^*_X (k)-\gamma_{1,X} (k)) \right \|_2 \leq \frac{C_1 a_n}{\sqrt  m} \, . 
$$
By condition \eqref{cond1} and the fact that $a_n=o(\sqrt n)$, this last quantity tends to $0$ as $n\rightarrow \infty$.

Let  now 
$$
\gamma_{2,X} (k)= \frac 1 n  \sum_{i=1}^{n-k}
 (H_m(X_i)-\bar H_m )(H_Y(X_{i+k}) - \pi)  \, .
$$
Then
\begin{equation}\label{rep2}
 \|\gamma_{1, X}(k)-\gamma_{2,X} (k)\|_2 \leq 
 \left \| \bar H_m-\frac 1 n \sum_{i=1}^n H_Y(X_i)
 \right \|_2 + \left \| \frac 1 n \sum_{i=1}^n H_Y(X_i) - \pi
 \right \|_2 \, .
\end{equation}
Proceeding as in \eqref{emp1} for the first term on the right-hand side of \eqref{rep2}, and noting that
$$
\left \| \frac 1 n \sum_{i=1}^n H_Y(X_i) - \pi
 \right \|_2^2 \leq  \frac{ 2 \sum_{k= 0}^{n-1}
\alpha_{1, \bf X}(k)}{n}\, ,
$$
we get  that there exist a positive constant $C_2$ such that
$$
\left \|\gamma_{1,X} (0)-\gamma_{2,X} (0) +2\sum_{k=1}^{a_n}  (\gamma_{1,X} (k)-\gamma_{2,X} (k)) \right \|_2 \leq \frac{C_2 a_n}{\sqrt  m} + \frac{C_2 a_n}{\sqrt  n} \, . 
$$
Again, this last quantity tends to $0$ as $n\rightarrow \infty$.

Finally, all the other replacements may be done in the same 
way, and lead to negligible contributions in ${\mathbb L}^2$.
Hence 
$$
\lim_{n \rightarrow \infty} \|V_n-V_n^*\|_2=0
$$
and the result follows from Step 1.

\subsection{Proof of Proposition \ref{BC}}
We start from \eqref{Hoef} again, but with the normalization 
$\sqrt{n / \log n}$ instead of $\sqrt n$.
Since \eqref{boundarycond}   holds,
$$
\limsup_{ n \rightarrow \infty}\frac 1 n \sum_{i=0}^n (i+1) \alpha_{1, \bf X}(i) < \infty \quad  \text{and} \quad 
\limsup_{ m \rightarrow \infty}\frac 1 m \sum_{j=0}^m (j+1) \alpha_{1, \bf Y}(j) < \infty \, .
$$
From Inequality \eqref{secondbound} and Condition \eqref{cond1}, 
we infer that
\begin{equation}\label{onebis}
\frac{\sqrt n}{n m \sqrt{\log n} } \sum_{i=1}^n \sum_{j=1}^m
f(X_i,Y_j) \quad \text{converges to 0 in ${\mathbb L}^2$ as $n \rightarrow \infty$.} 
\end{equation}

Let us prove Item 1. From \eqref{Hoef} and \eqref{onebis},
\begin{multline}\label{dev}
\limsup_{n \rightarrow \infty}
  {\mathbb P} \left ( \sqrt{n /\log n} (U_n-\pi)>x \right)
   \leq 
  \limsup_{n \rightarrow \infty}
  {\mathbb P} \left (  \sum_{i=1}^n (H_Y(X_i)- \pi) > x \sqrt{n \log n} /2 \right) \\
  + \limsup_{n \rightarrow \infty}
  {\mathbb P} \left ( \frac n m \sum_{j=1}^m (G_X(Y_i)- \pi) > x \sqrt{n \log n} /2 \right) \, .
\end{multline}
Since \eqref{boundarycond} holds, it follows from Proposition 5.1
in \cite{DM} that: for any $r \in (2,4)$, there exists a positive constant $c_{1,r}$ such that, for any $x>0$
$$
\limsup_{n \rightarrow \infty}
  {\mathbb P} \left (  \sum_{i=1}^n (H_Y(X_i)- \pi) > x \sqrt{n \log n} /2 \right) \leq  \frac{c_{1,r}}{x^r} \, .
$$
Since \eqref{cond1} holds, the same is true for the second term
on the right-hand side of \eqref{dev} (for a positive 
constant  $c_{2,r}$), and Item 1 follows.

Let us prove Item 2. We use a result by Gou\"ezel \cite{G} about intermittent maps. Let $\theta_{1/2}$ be the map described in Subsection \ref{Sec:IM} of the present paper.
As explained in Subsection \ref{Sec:IM}, there exists an unique absolutely continuous $\theta_{1/2}$-invariant measure $\nu_{1/2}$ on 
$[0,1]$, and (see the comments on the case $\alpha=1/2$ in Section 1.3 of Gou\"ezel's paper), on  the probability space $([0,1], \nu_{1/2})$,
$$
\frac {1}{\sqrt n} \sum_{k=1}^n\left ( H_Y(\theta_{1/2}^k) -\nu_{1/2} ( H_Y)
\right )$$
converges in distribution to a non-degenerate normal distribution provided $H_Y(0)\neq \nu(H_Y)$ (which is true
for instance if $Y$ has distribution $\nu_{1/2}$, since in that 
case $H_Y(0)=1$ and $\nu(H_Y)=1/2$).
As fully explained in Subsection  \ref{Sec:IM}, there exists 
a stationary  Markov chain $(X_i)_{i \in {\mathbb Z}}$ such
that the distribution of $(X_1, \dots, X_n)$ is 
the same as $(\theta_{1/2}^n, \ldots, \theta_{1/2})$ on the probability  space 
$([0,1], \nu_{1/2})$. Take $(Y_i)_{i \in {\mathbb Z}}$ an independent
copy of the chain $(X_i)_{i \in {\mathbb Z}}$.  For such chains, 
both 
$$
\frac{1}{ \sqrt n} \sum_{i=1}^n (H_Y(X_i)- \pi)  \quad \text{and}
\quad 
\frac {1}{\sqrt m} \sum_{j=1}^m (G_X(Y_j)- \pi)
$$
converge to the same non-degenerate normal distribution
(note that $\pi=1/2$ in that case). Moreover, as recalled in  Subsection  \ref{Sec:IM},
there exist two positive constants 
$A,B$ such that 
$$
\frac{A}{k+1} \leq \alpha_{2, \bf X}(k)=\alpha_{2, \bf Y}(k) \leq \frac{B}{k+1} \, .
$$
This completes the proof of Item 2.

\subsection{Proof of Proposition \ref{LR}}
We start from \eqref{Hoef} again. At the end of the proof, we shall prove that if \eqref{LRcond} and \eqref{cond2} are satisfied, then
\eqref{one} holds. 

Let us now prove  Item 1. From \eqref{Hoef} and \eqref{one},
\begin{multline}\label{devbis}
\limsup_{n \rightarrow \infty}
  {\mathbb P} \left ( \sqrt{n } (U_n-\pi)>x \right)
   \leq 
  \limsup_{n \rightarrow \infty}
  {\mathbb P} \left (  \sum_{i=1}^n (H_Y(X_i)- \pi) > x \sqrt{n} /2 \right) \\
  + \limsup_{n \rightarrow \infty}
  {\mathbb P} \left ( \frac n m \sum_{j=1}^m (G_X(Y_i)- \pi) > x \sqrt{n } /2 \right) \, .
\end{multline}
Since the first part of  \eqref{LRcond} is satisfied, it follows from
the central limit theorem for  stationary and 2-ergodic $\alpha$-dependent sequences that
\begin{equation}\label{Gaussdev}
\lim_{n \rightarrow \infty}
  {\mathbb P} \left (  \sum_{i=1}^n (H_Y(X_i)- \pi) > x \sqrt{n } /2 \right) ={\mathbb P}(Z>x/2) \, ,
\end{equation}
where $Z$ is a Gaussian random variable. From Condition \eqref{LRcond}  and Remark 3.3 
in \cite{DM}, we infer  that  if \eqref{cond2} holds, then there exists a positive constant $\kappa_p$ such that, for any $x>0$,
\begin{equation}\label{Stabledev}
\limsup_{n \rightarrow \infty}
{\mathbb P} \left ( \frac n m \sum_{j=1}^m (G_X(Y_i)- \pi) > x \sqrt{n } /2 \right) \leq \frac{\kappa_p}{x^p} \, .
\end{equation}
Item 1 follows from \eqref{devbis}, \eqref{Gaussdev} and 
\eqref{Stabledev}.

Let us prove Item 2.
For the sequence $(X_i)_{i \in {\mathbb Z}}$, take a sequence 
of iid random variables uniformly distributed over $[0,1]$. 
To find the sequence $(Y_i)_{i \in {\mathbb Z}}$  we use again  a result by Gou\"ezel \cite{G} about intermittent maps. Let $\theta_{1/p}$ be the map described in Subsection  \ref{Sec:IM} 
of the present paper. As explained is Subsection \ref{Sec:IM}, there exists an unique absolutely continuous $\theta_{1/p}$-invariant measure $\nu_{1/p}$ on 
$[0,1]$, and (see Theorem 1.3 in \cite{G}), on  the probability space $([0,1], \nu_{1/p})$,
$$
\frac {1}{m^{1/p}} \sum_{k=1}^m\left ( G_X(\theta_{1/p}^k) -\nu_{1/p} ( G_X)
\right )$$
converges in distribution to a non-degenerate stable distribution of order $p$ provided $G_X(0)\neq \nu(G_X)$ (this condition 
is satisfied here, because $G_X(0)=0$ and $\nu_{1/p}(G_X)>0$).
As fully explained in Subsection  \ref{Sec:IM}, there exists 
a stationary  Markov chain $(Y_i)_{i \in {\mathbb Z}}$ such
that the distribution of $(Y_1, \dots, Y_n)$ is 
the same as $(\theta_{1/p}^n, \ldots, \theta_{1/p})$ on the probability  space 
$([0,1], \nu_{1/p})$.   For such  a chain, 
$$
\frac {1}{m^{1/p}} \sum_{j=1}^m (G_X(Y_j)- \pi)
$$
converges in distribution  to  a non-degenerate stable distribution
of  order $p$. Moreover, as recalled in  Subsection  \ref{Sec:IM}, there exist two positive constants 
$A,B$ such that 
$$
\frac{A}{(k+1)^{p-1}} \leq \alpha_{1, \bf Y}(k) \leq \frac{B}{(k+1)^{p-1}} \, .
$$
This completes the proof of Item 2.

It remains to prove \eqref{one}. From Proposition \ref{degen} we get that 
\begin{equation}\label{equiv}
n{\mathbb E}\left (\left( \frac{1}{nm} \sum_{i=1}^n \sum_{j=1}^m
f(X_i,Y_j) \right )^2 \right ) \leq 
\sum_{i=0}^{\infty} \left (\frac{16}{m}\sum_{j=0}^{m-1} \min \{ \alpha_{1, \bf X}(i), \alpha_{1, \bf Y}(j) \}
\right )
  \, .
\end{equation}
Note that 
\begin{equation}\label{borne}
\frac{16}{m}\sum_{j=0}^{m-1} \min \{ \alpha_{1, \bf X}(i), \alpha_{1, \bf Y}(j) \} \leq 16 \alpha_{1, \bf X}(i) \, ,
\end{equation}
and that, for any positive integer $i$, since 
$\alpha_{1, \bf Y}(j) \rightarrow 0$ as $j \rightarrow \infty$, 
\begin{equation}\label{lim}
\lim_{m \rightarrow \infty} \frac{16}{m}\sum_{j=0}^{m-1} \min \{ \alpha_{1, \bf X}(i), \alpha_{1, \bf Y}(j) \}  =0 \, ,
\end{equation}
by Cesaro's lemma. Since $\sum_{i>0} \alpha_{1, \bf X}(i) < \infty$, it follows from \eqref{equiv}, \eqref{borne}, 
\eqref{lim}, and the dominated convergence theorem 
that 
$$
\lim_{n \rightarrow \infty} 
n {\mathbb E}\left ( \left ( 
\frac{1}{nm}  \sum_{i=1}^n \sum_{j=1}^m
f(X_i,Y_j)\right )^2 \right ) =0 \, ,
$$
proving \eqref{one}.

\section{Other related results}

\subsection{Testing the stochastic domination by a known distribution}
\setcounter{equation}{0}
We briefly describes here a procedure for testing the 
(weak) domination by a known distribution $\mu$. 
Let $H(t)=\mu((t, \infty))$, and let $(X_i)_{i \in {\mathbb Z}}$ 
be a stationary 2-ergodic sequence of real-valued random variables. 
Let $\pi={\mathbb E}(H(X_0))$, and define  
$$
\bar{H}_n= \frac 1 n \sum_{i=1}^n H(X_i) \, .
$$
If 
\begin{equation}\label{alphasimple}
\sum_{k=1}^\infty \alpha_{1, \bf X}(k) < \infty \, ,
\end{equation}
then the random variables
$
\sqrt n (\bar{H}_n - \pi)
$
converge in distribution to ${\mathcal N}(0, V_1)$, where
\begin{equation}\label{Var1}
V_1= \mathrm{Var}(H(X_0)) + 2 \sum_{k=1}^\infty 
\mathrm{Cov}(H(X_0), H(X_k)) \, .
\end{equation}

Now, we want to test $H_0: \pi=1/2$ (no relation of weak domination between $\mu$ and ${\mathbb P}_X$)
against one of the possible alternatives:
$H_1: \pi \neq 1/2$ (one of the two distributions 
weakly dominates the other one), $H_1: \pi>1/2$ ($\mu$ weakly 
dominates ${\mathbb P}_X$), $H_1: \pi<1/2$ (${\mathbb P}_X$
weakly dominates $\mu$). 

Again, one cannot use directly the statistics 
$\sqrt{n}(\bar{H}_n-1/2)$ to test $H_0: \pi=1/2$ versus one of the possibles alternatives, because under $H_0$ the asymptotic distribution of $\sqrt{n}(\bar{H}_n-1/2)$ depends of the unknown 
quantity $V_1$ defined in \eqref{Var1}. 

To estimate $V_1$, we take 
$$
V_{1,n}=  \hat \gamma (0) + 2\sum_{k=1}^{a_n} \hat \gamma (k)  
$$
for some sequence $a_n=o(\sqrt n)$, where
$$
\hat \gamma (k)= \frac 1 n  \sum_{i=1}^{n-k}
 (H(X_i)-\bar{H}_n )(H(X_{i+k}) - \bar{H}_n)  \, . 
$$
Now, as in Proposition \ref{varest}, if \eqref{alphasimple} holds 
for $\alpha_{2, \bf X}(k)$ instead of $\alpha_{1, \bf X}(k)$, then 
$v_n$ converges in ${\mathbb L}^2$ to $v$. 
It follows that, under $H_0: \pi=1/2$, if $v>0$, the random variables
$$
T'_n= \frac{\sqrt{n}(\bar{H}_n - 1/2)}{\sqrt{\max \{V_{1,n}, 0\}}}
$$
converges in distribution to ${\mathcal N}(0,1)$. 

\subsection{The case of two adjacent samples from 
the same time series}\label{adj}

Let $(X_i)_{i \in {\mathbb Z}}$ 
be a stationary sequence of real-valued random variables. Let $m=m_n$ be 
a sequence of positive integers such that $m_n \rightarrow \infty$ as $n \rightarrow \infty$. 
We consider here the   U-statistic
$$
U_n= U_{n,m}= \frac{1}{n m } 
\sum_{i=1}^n \sum_{j=n+1}^{n+m}  {\bf 1}_{X_i < Y_j}  \, , 
$$
where $Y_j=f(X_j)$, $f$ being a monotonic function. 
This statistic can be used to test if there is a change 
in the time series after the time $n$, and more precisely 
if $Y_{n+1}$ weakly dominates $X_1$. Let 
$\pi={\mathbb P}(X_1 < f(X_1^*))$, where 
$X_1^*$ is an independent copy
of $X_1$.
Let also 
$H_Y(t)={\mathbb P}(f(X_1)>t)$ and 
$G_X(t)={\mathbb P}(X_1 <t)$. 

Finding the asymptotic behavior of $U_n$ is a slightly more difficult problem than for two independent samples, 
because we cannot use the independence to control
the degenerate U-statistic as in Proposition 
\ref{degen}. However, it can be done 
by using a more restrictive coefficient than 
${\alpha}_{2, {\bf X}}$. 

\begin{Definition}\label{defbeta}
 Let $(X_i)_{i \in {\mathbb Z}}$ be a strictly stationary sequence
 of real-valued random variables, 
 and let ${\mathcal F}_0=\sigma (X_k, k \leq 0)$. 
 Let $P$ be the law of $X_0$ and $P_{(X_i, X_j)}$ be the law of $(X_i, X_j)$. Let
 $P_{X_k| {\mathcal
F}_0}$ be the conditional distribution of  $X_k$ given ${\mathcal F}_0$, and  let $P_{(X_i,
X_j)|{\mathcal F}_0}$ be the conditional distribution of $(X_i,X_j)$ given ${\mathcal
F}_0$. Define  the
  random variables
\begin{eqnarray*}
b (k) &=& \sup_{x \in {\mathbb R}}
\left |{P}_{X_{k}|{\mathcal F}_0}(f_x) \right |\, , \\
b ( i,j)&=& \sup_{(x,y) \in {\mathbb R}^2} \left |P_{(X_i, X_j)|{\mathcal F}_0}\left (f_x\otimes f_y \right )-P_{(X_i, X_j)}\left (f_x\otimes f_y\right ) \right |\, ,
\end{eqnarray*}
where the function $f_x$ has been defined in Definition \ref{defalpha}.
Define now the coefficients
\[
 \beta_{1,{\bf X}}(k)= {\mathbb E}(b( k))\,  \text{ and } \,   \beta_{2,{\bf X}}(k) =\max \left \{\tilde \beta_{1, {\bf X}}(k), \sup_{i>j\geq k} {\mathbb E}((b( i,j))) \right \}\, .
\]
\end{Definition}

\smallskip

Of course, by definition of the coefficients, we always have that
$\alpha_{2, {\bf X}}(k) \leq \beta_{2, {\bf X}}(k)$.
The coefficients $\beta_{2, {\bf X}}$ are weaker than the usual $\beta$-mixing coefficients of
 $(X_i)_{i \in {\mathbb Z}}$. Many examples of non-mixing process 
for which $ \beta_{2,{\bf X}}$ can be computed are given in \cite{DP2}. For  all the examples that we mention in the present paper, the coefficient $\beta_{2, {\bf X}}(k)$ is of the same order
than the coefficient $\alpha_{2, {\bf X}}(k)$ (see the paper 
 \cite{DDT} for the example of the intermittent maps  described in  Subsection \ref{Sec:IM}).

 \begin{Theorem}\label{MannThbis}
Let $(X_i)_{i \in {\mathbb Z}}$ 
be a stationary sequence
of real-valued random variables. Assume that 
\begin{equation}\label{betacond}
\sum_{k=1}^\infty \beta_{2,\bf{X}}(k) < \infty  \, .
\end{equation} 
Let $(m_n)$ be a sequence of positive integers 
satisfying \eqref{cond1}.
Then 
$
  \sqrt{n} (U_n- \pi) 
$
converge in distribution to ${\mathcal N}(0,V)$, with
\begin{multline}\label{Varbis}
V=\mathrm{Var}(H_Y(X_0)) + 2 \sum_{k=1}^\infty 
\mathrm{Cov}(H_Y(X_0), H_Y(X_k)) \\ + 
c \left ( \mathrm{Var}(G_X(f(X_0))) + 2 \sum_{k=1}^\infty 
\mathrm{Cov}(G_X(f(X_0)), G_X(f(X_k)))\right ) \, .
\end{multline}
\end{Theorem}


For the estimation of $V$, one can prove the following 
analogue of Proposition \ref{varest}. 

\begin{Proposition}\label{varest2}
Let $(X_i)_{i \in {\mathbb Z}}$ 
be an independent stationary sequence, and let $Y_j=f(X_j)$ 
where $f$ is a monotonic function.
Let $(m_n)$ be a sequence of positive integers satisfying
\eqref{cond1}.  Assume that \eqref{betacond} is satisfied.
Let 
$$
H_m(t)= \frac 1 m \sum_{j=n+1}^{n+m} {\bf 1}_{Y_j>t} \quad \text{and}
\quad 
G_n(t)= \frac 1n \sum_{i=1}^n {\bf 1}_{X_i<t} \, , 
$$
and 
$$
 \hat \gamma_X (k)= \frac 1 n  \sum_{i=1}^{n-k}
 (H_m(X_i)-\bar{H}_m )(H_m(X_{i+k}) - \bar{H}_m) \, , \,
  \hat \gamma_Y (\ell)= \frac 1 m  \sum_{j=n+1}^{n+m-\ell}
 (G_n(Y_j)-\bar{G}_n )(G_n(Y_{j+\ell}) - \bar{G}_n) \, ,
$$
where $n \bar{H}_m=  \sum_{i=1}^{n}
 H_m(X_i)$ and $m \bar{G}_n=  \sum_{j=n+1}^{n+m}
 G_n(Y_j)$. 
Let also $(a_n)$ and $(b_n)$ be two sequences of positive
 integers tending to infinity as $n$ tends to infinity, such that $a_n=o(\sqrt n/\log(n))$ and $b_n=o(\sqrt {m_n}/\log(m_n))$. 
 Then 
 $$
 V_n=  \hat \gamma_X (0) +2 \sum_{k=1}^{a_n} \hat \gamma_X (k)  
 +\frac{n}{m_n} \left ( 
 \hat \gamma_Y (0) +2 \sum_{\ell=1}^{b_n} \hat \gamma_Y (\ell)   \right )
 $$
 converges in ${\mathbb L}^2$ to the quantity $V$ defined
 in \eqref{Varbis}.
\end{Proposition}

Combining Theorem \ref{MannThbis} and Proposition \ref{varest2}, 
we obtain that, under $H_0: \pi=1/2$, if $V>0$, the random variables
\begin{equation}\label{Tnbis}
T_n= \frac{\sqrt{n}(U_n - 1/2)}{\sqrt{\max \{V_n, 0\}}}
\end{equation}
converge in distribution to ${\mathcal N}(0,1)$. 

\medskip

The proof of Proposition \ref{varest2} follows the same lines as that of Proposition \ref{varest}, so we do not give the details here. The only main change concerns
Step 2, Inequalities \eqref{rep1} and \eqref{emp1}, where we used
the independence between the samples. Here instead, 
one can use the fact that 
$$
\left \|\sup_{x \in {\mathbb R}}|H_n(x)-H(x)|
\right  \|_2^2\leq 
C \left (\log(m)
\right)^2\frac{\sum_{k = 0}^{m-1} \beta_{1, {\bf X}}(k)}
{m} \, .
$$
This can be easily proved by following the proof of 
Proposition 7.1 in Rio \cite{Ri} up to (7.8) and by using the computations   in the proof of Theorem 2.1
in \cite{D}. 

\smallskip

The proof of Theorem \ref{MannThbis} follows also 
the line of that of Theorem \ref{MannTh}; the only non-trivial change concerns Proposition \ref{degen}.
Here, we have to deal with 
$$
  L_n= \frac{1}{nm} \sum_{i=1}^n \sum_{j=n+1}^{n+m}
f(X_i,Y_j) 
$$
where $f$ is defined in \eqref{degenU}. We want to  prove 
that $n{\mathbb E}(L_n^2)$ tends to 0 as $n$ tends 
to infinity. As in the proof of Proposition \ref{degen}, 
we start from the elementary inequality
\begin{equation} \label{elem2}
{\mathbb E}(L_n^2)
\leq \frac{4}{n^2 m^2} \sum_{i=1}^n \sum_{j=n+1}^{n+m}
\sum_{k=i}^n \sum_{\ell =j}^{n+m}
{\mathbb E}( f(X_i, Y_j) f(X_k, Y_\ell)) \, .
\end{equation}
Let $(X_i^*)_{ i \in {\mathbb Z}}$ be an
independent copy of $(X_i^*)_{ i \in {\mathbb Z}}$
and $Y_i^*=f(X_i^*)$, and write
\begin{multline}\label{trivdev}
{\mathbb E}( f(X_i, Y_j) f(X_k, Y_\ell))=
{\mathbb E}( f(X_i, Y_j) f(X_k, Y_\ell))-
{\mathbb E}( f(X_i, Y_j^*) f(X_k, Y_\ell^*)) 
\\
+{\mathbb E}( f(X_i, Y_j^*) f(X_k, Y_\ell^*)) \, .
\end{multline}
Now, using the independence, the second term
on the right-hand side of \eqref{trivdev} can be handled 
exactly as in Proposition \ref{degen}. It remains to 
deal with the first term in the right-hand side of 
\eqref{trivdev}. Note that
\begin{multline}
|{\mathbb E}( f(X_i, Y_j) f(X_k, Y_\ell))-
{\mathbb E}( f(X_i, Y_j^*) f(X_k, Y_\ell^*)) | \\
\leq 
\int |{\mathbb E}( f(x_i, Y_j) f(x_k, Y_\ell)| X_i=x_i, X_k=x_k) -
{\mathbb E}( f(x_i, Y_j) f(x_k, Y_\ell)) | P_{(X_i, X_k)} (dx_i, 
dx_k)
\end{multline}
For any 
fixed $x$, the function $y \rightarrow h_x(y)=f(x,y)$ is a bounded variation function whose variation is bounded
by 2, and ${\mathbb E}(h_x(X_i))=0$. 
Proceeding as in Lemma 1 
of \cite{DP1}, 
one can then prove that 
\begin{multline}
\int |{\mathbb E}( f(x_i, Y_j) f(x_k, Y_\ell)| X_i=x_i, X_k=x_k) -
{\mathbb E}( f(x_i, Y_j) f(x_k, Y_\ell)) | P_{(X_i, X_k)} (dx_i, 
dx_k) \\
\leq 4 \beta_{2, {\bf X}}(j-k) \, .
\end{multline}
Considering separately the two cases $j-k> n^{1/4}$ or 
$j-k \leq n^{1/4}$, we get the upper bound
\begin{multline}
\frac{1}{n^2 m^2} \sum_{i=1}^n \sum_{j=n+1}^{n+m}
\sum_{k=i}^n \sum_{\ell =j}^{n+m}
|{\mathbb E}( f(X_i, Y_j) f(X_k, Y_\ell)-
{\mathbb E}( f(X_i, Y_j^*) f(X_k, Y_\ell^*))| \\ \leq 
\frac{4}{\sqrt n m } + \frac{4}{m} \sum_{k > n^{1/4}}
 \beta_{2, {\bf X}}(k)
\end{multline}
from which we easily derive that $n{\mathbb E}(L_n^2)$ tends to 0 as $n$ tends 
to infinity (since \eqref{cond1} and \eqref{betacond} are satisfied).

\section{Simulations}\label{Sec:Simu}

\setcounter{equation}{0}

\subsection{Example 1:  functions of an auto-regressive process}\label{secar}
In this section, 
we first simulate $(Z_1, \ldots, Z_n)$,  according to the simple 
AR(1) equation
$$
Z_{k+1}=\frac12 \left ( Z_k + \varepsilon_{k+1} \right ) \, ,
$$
where $Z_1$ is uniformly distributed over $[0,1]$,
and $(\varepsilon_i)_{i \geq 2}$ is a sequence of 
iid random variables with distribution ${\mathcal B}(1/2)$,
independent of $Z_1$.

One can check that the transition Kernel of the  chain  $(Z_i)_{i \geq 1}$ is
$$
 K(f)(x)= \frac 1 2 \left (f\left (\frac x 2\right )+ f\left (\frac{x+1}{2} \right ) \right ) \, ,
$$
and that the uniform distribution on $[0,1]$ is the unique
invariant distribution by $K$. Hence, the chain $(Z_i)_{i \geq 1}$
is strictly stationary. 

It is well known that the chain $(Z_i)_{i \geq 1}$ is not $\alpha$-mixing 
in the sense of Rosenblatt \cite{R} (see for instance \cite{Br}). 
However, one can  prove that the coefficients $\alpha_{2,\bf X}$
of the chain $(Z_i)_{i \geq 1}$ are such that 
$$
   \alpha_{2, {\bf Z}}(k) \leq 2^{-k} 
$$
(and the same upper bound is true for $\beta_{2, {\bf Z}}(k)$, see for instance 
Section 6.1 in \cite{DP2}). 

\smallskip

Let now $Q_{\mu, \sigma^2}$ be the inverse of the 
cumulative distribution function of the law 
${\mathcal N}(\mu, \sigma^2)$. Let then
$$
X_i=Q_{\mu_1, \sigma_1^2}(Z_i) \, .
$$
The sequence $(X_i)_{i \leq 1}$ is also a stationary Markov chain (as 
an invertible function of a stationary Markov chain), and one can 
easily check that 
$ \alpha_{2,{\bf X}}(k)= \alpha_{2,{\bf Z}}(k)$. By construction, 
$X_i$ is ${\mathcal N}(\mu, \sigma^2)$-distributed, 
but the sequence $(X_i)_{i \geq 1}$ is not a Gaussian 
process (otherwise it would be mixing in the sense of
Rosenblatt). 

\smallskip

Next, we simulate iid random variables 
$(Y_1, \ldots, Y_m)$   with
common distribution ${\mathcal N}(\mu_2, \sigma_2^2)$.

\smallskip

Note that, if $\mu_1=\mu_2$, the hypothesis
$H_0: \pi=1/2$  is satisfied and the random variables $T_n$ defined in \eqref{Tn} converge in  
distribution to ${\mathcal N}(0,1)$. 

\smallskip

We shall  study the behavior 
of $T_n$  for  $\sigma_1^2=4$, $\sigma_2^2=1$ and different choices of $\mu_1, \mu_2, a_n$ and $b_n$. As explained in 
Section \ref{sec2}, the quantities $a_n$ and $b_n$ may be chosen by analyzing the 
graph of the estimated auto-covariances $\hat \gamma_X(k)$ and $\hat \gamma_Y(\ell)$. 
For this examples, these graphs suggest a choice of $a_n=3$ or $4$ and $b_n=0$
(see Figure 1). 

\smallskip

We shall  compute $T_n$ for different choices of $(n,m)$ (from $(150,100)$ to 
$(750,500)$ with $n/m=1.5$). We estimate the three quantities 
$\mathrm {Var}(T_n)$, ${\mathbb P}(T_n>1.645)$  and 
 ${\mathbb P}(|T_n|>1.96)$  via a classical Monte-Carlo procedure, by averaging
 over $N=2000$ independent trials.
 
 \smallskip
 
  If $\mu_1= \mu_2$,  we say that
${\mathbb P}(T_n>1.645)$ is level 1 (the level of the 
one-sided test)  and ${\mathbb P}(|T_n|>1.96)$ is level 2
(the level of the two-sided test). If 
  $a_n, b_n$ are well chosen, the estimated variance
 should be close to 1 and the estimated levels 1 and 2 should be close to $0.05$. 
 
 \smallskip
 
 If $\mu_1\neq  \mu_2$, we say that
${\mathbb P}(T_n>1.645)$ is power 1 (the power of the 
one-sided test)  and ${\mathbb P}(|T_n|>1.96)$ is power
2 (the power of the two-sided test).

 \medskip

$\bullet$ Case $\mu_1=\mu_2=0$ and  $a_n=0, b_n=0$ (no correction):

\medskip

\begin{center}
\begin{tabular}{|c|c|c|c|c|c|}
\hline
\quad \quad \quad \quad \quad \quad $n ; m$  
&    150 ; 100  & 300 ; 200 & 450 ; 300 & 600 ; 400 & 750 ; 500 \\
\hline 
Estimated variance &2.13  & 2.199  &  2.118 & 2.112 &  2.273 \\
\hline
Estimated level 1 &0.129 &0.145 & 0.121& 0.132 & 0.14\\
\hline
Estimated level 2 &0.177 &0.194 & 0.18 & 0.18& 0.191 \\
\hline
\end{tabular}
\end{center}

\medskip

Here, since $a_n=0$, we do not estimate any of the 
covariance terms $\gamma_X(\ell)$; hence, we compute 
$T_n$ as if both series $(X_i)_{i \geq 1}$ and $(Y_j)_{j \geq 1}$ were uncorrelated. 
The result is that the estimated variance is too large (around 2.1 whatever $(n,m)$)
and so are the estimated levels 1 and 2 (around 0.13 for level 1 and 0.18 for level 2, 
whatever $(n,m)$). This means that the-one sided or two-sided tests build with $T_n$ 
will reject the null hypothesis too often.

\medskip

\begin{figure}[htb]
\centerline{
\includegraphics[width=11.5 cm, height=6.5 cm]{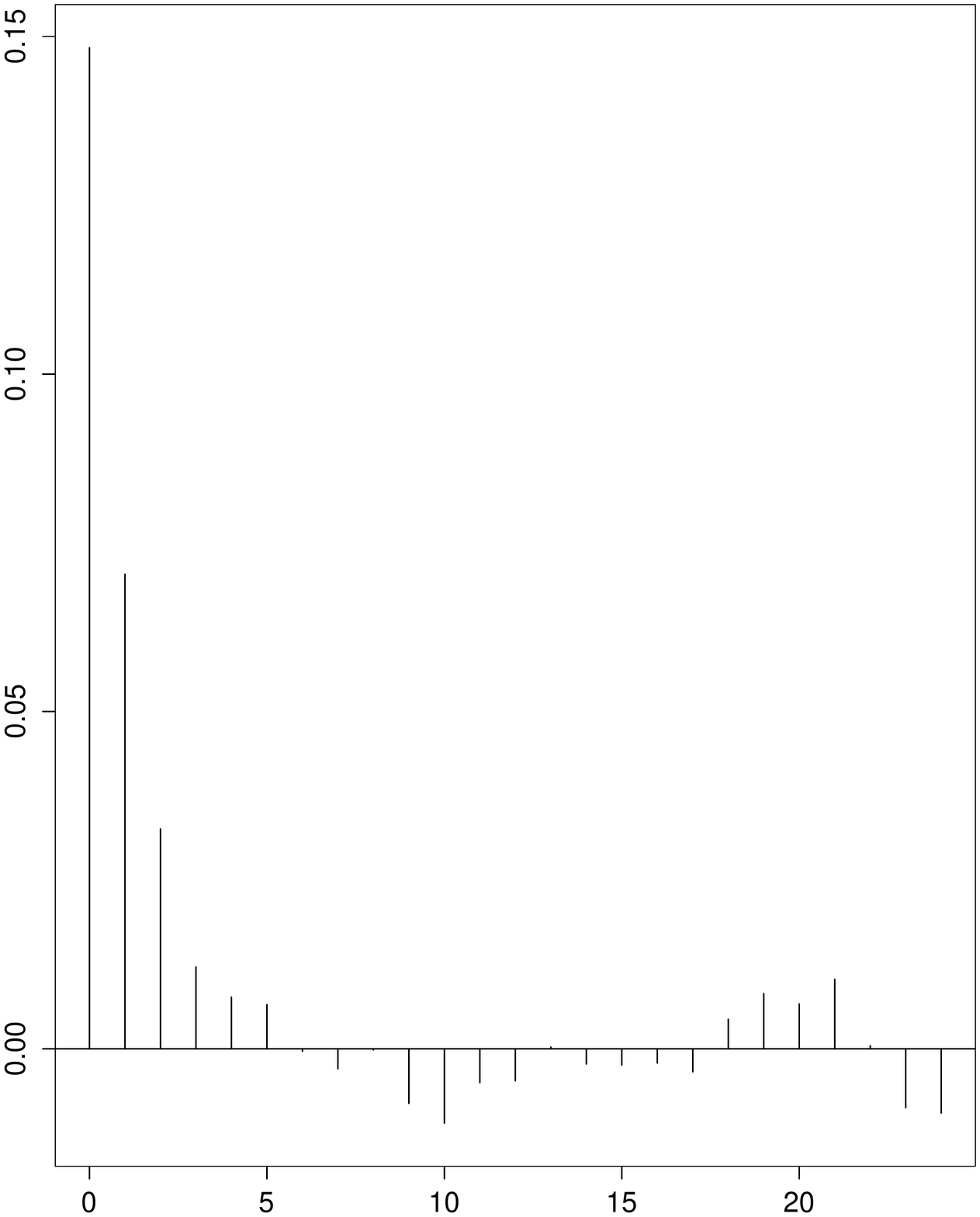}}
\centerline{
\includegraphics[width=11.5 cm, height=6.5 cm]{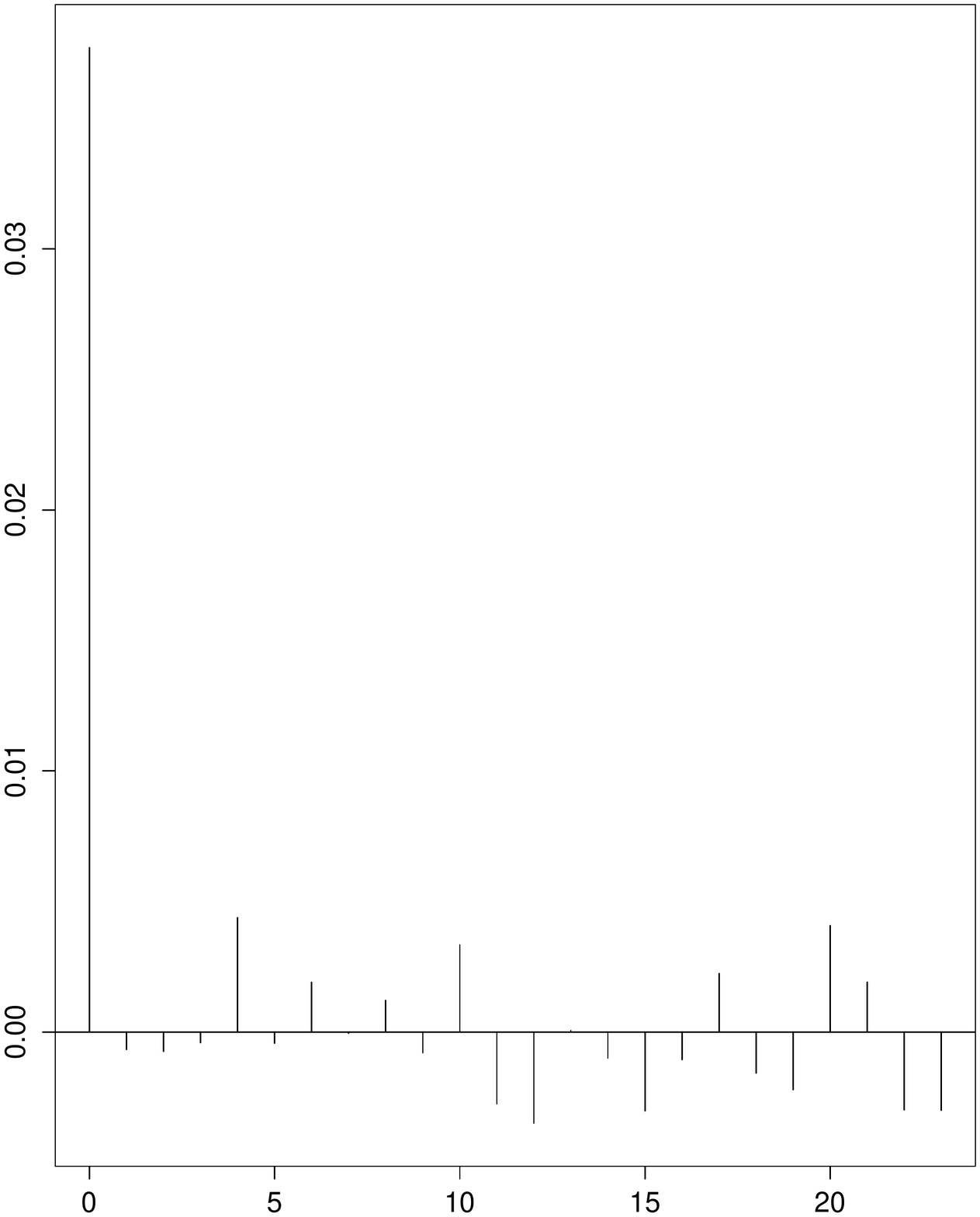}}
\caption{Estimation of the $\gamma_X(k)$'s (top) and 
$\gamma_Y(\ell)$'s (bottom)
for Example 1, 
with   $\mu_1=\mu_2=0$, $n=300$  and $m=200$.}
\label{fig:covar}
\end{figure}

\medskip

$\bullet$ Case $\mu_1=\mu_2=0$ and  $a_n=3, b_n=0$:

\medskip

\begin{center}
\begin{tabular}{|c|c|c|c|c|c|}
\hline
\quad \quad \quad \quad \quad \quad $n ; m$  
&    150 ; 100  & 300 ; 200 & 450 ; 300 & 600 ; 400 & 750 ; 500 \\
\hline 
Estimated variance &1.227  & 1.148  &  1.172 & 1.119 &  1.018 \\
\hline
Estimated level 1 &0.071 &0.061 & 0.065& 0.06 & 0.046\\
\hline
Estimated level 2 &0.074 &0.069 & 0.071 & 0.066& 0.054 \\
\hline
\end{tabular}
\end{center}

\medskip

As suggested by the graphs of the estimated auto-covariances (see Figure 1), the 
choice   $a_n=3, b_n=0$ should give a reasonable estimation of $V$. This is indeed
the case: the estimated variance is now around $1.1$, the estimated level 1 is around
$0.06$ and the estimated level 2 around $0.07$.

\medskip

$\bullet$  Case $\mu_1=\mu_2=0$ and  $a_n=4, b_n=0$:

\medskip 

\begin{center}
\begin{tabular}{|c|c|c|c|c|c|}
\hline
\quad \quad \quad \quad \quad \quad $n ; m$  
&    150 ; 100  & 300 ; 200 & 450 ; 300 & 600 ; 400 & 750 ; 500 \\
\hline 
Estimated variance &1.151  & 1.11  &  1.121 & 1.03 &  0.987 \\
\hline
Estimated level 1 &0.06 &0.05 & 0.058 & 0.049 & 0.042\\
\hline
Estimated level 2 &0.062 &0.067 & 0.06 & 0.057& 0.046 \\
\hline
\end{tabular}
\end{center}

\medskip

Here, we see that the choice $a_n=4, b_n=0$ works well
also, and seems even slightly better than the choice $a_n=3, b_n=0$.  

\medskip

$\bullet$ Case $\mu_1=0$, $\mu_2=0.2$ and  $a_n=4, b_n=0$ :

In this example, $H_0$ is not satisfied, and $\pi>1/2$. However, $\mu_2$ is close to 
$\mu_1$, and since $\sigma_1=2$ and $\sigma_2=1$, $Y_1$ does not stochastically 
dominate $X_1$ (there is only stochastic domination in a weak sense). We shall
see how the two tests are able to see this weak domination, by giving some estimation 
of the powers.

\medskip

\begin{center}
\begin{tabular}{|c|c|c|c|c|c|}
\hline
\quad \quad \quad \quad \quad \quad $n ; m$  
&    150 ; 100  & 300 ; 200 & 450 ; 300 & 600 ; 400 & 750 ; 500 \\
\hline 
Estimated variance &1.229  & 1.136  &  1.107 & 1.066 &  1.059 \\
\hline
Estimated power 1 &0.185 &0.276 & 0.352 & 0.4 & 0.476\\
\hline
Estimated power 2 &0.13 &0.187 & 0.25  & 0.291& 0.344 \\
\hline
\end{tabular}
\end{center}

\medskip

As one can see, the powers of the one sided and two-sided tests are always greater than $0.05$, as 
expected. Still as expected, these powers increase with the size of the samples. 
For $n=750, m=500$, the power of the one-sided test is around $0.47$
and that of the two-sided test is around $0.34$.

\subsection{Example 2: Intermittent maps} \label{Sec:IM}

For  $\gamma$ in $]0, 1[$, we consider the intermittent map $\theta_\gamma$
from $[0, 1]$ to
$[0, 1]$, introduced   by Liverani, Saussol and Vaienti \cite{LSV}:
\begin{equation}\label{LSV}
   \theta_\gamma(x)=
  \begin{cases}
  x(1+ 2^\gamma x^\gamma) \quad  \text{ if $x \in [0, 1/2[$}\\
  2x-1 \quad \quad \quad \ \  \text{if $x \in [1/2, 1]$.}
  \end{cases}
\end{equation}
It follows from  \cite{LSV} that
there exists a
unique absolutely continuous $\theta_\gamma$-invariant probability measure $\nu_\gamma$,  
with density $h_\gamma$. 

Let us briefly describe the Markov chain associated with $\theta_\gamma$, and its properties.
Let first $K_\gamma$ be the
 Perron-Frobenius operator of $\theta_\gamma$ with respect to $\nu_\gamma$, defined as follows:
for any  functions $u, v$ in ${\mathbb L}^2([0,1], \nu_\gamma)$
\begin{equation}\label{Perron}
\nu_\gamma (u \cdot  v\circ \theta_\gamma)=\nu_\gamma (K_\gamma (u) \cdot  v) \, .
\end{equation}
The relation \eqref{Perron} states  that $K_\gamma$ is  the adjoint operator of the 
isometry $U: u \mapsto u\circ \theta_\gamma$
acting on ${\mathbb L}^2([0,1], \nu_\gamma)$. It is easy to see that the operator $K_\gamma$ is a transition kernel, 
and that $\nu_\gamma$ is invariant by $K_\gamma$. 
Let now 
$(\xi_i)_{i\geq 1}$
be a stationary Markov chain with invariant measure $\nu_\gamma$
and transition kernel $K_\gamma$. It is well known  that on the probability space 
$([0, 1], \nu_\gamma)$, the
random vector $(\theta_\gamma, \theta_\gamma^2, \ldots , \theta_\gamma^n)$ is distributed as
$(\xi_n,\xi_{n-1}, \ldots, \xi_1)$. Now it is proved in \cite{DGM} that there exist
two positive constants $A,B$ such that
$$
\frac {A}{(n+1)^{(1-\gamma)/\gamma}} \leq \alpha_{2, \bf \xi}(n) 
\leq \frac{B}{(n+1)^{(1-\gamma)/\gamma}}\, 
$$
(this is also true for the coefficient $\beta_{2, \bf \xi}(n)$, see \cite{DDT}). 
Moreover, the chain $(\xi_i)_{i\geq 1}$ is not $\alpha$-mixing in the 
sense of Rosenblatt \cite{R}. 

\smallskip

Let $X_i(x)=\theta_{\gamma_1}^i(x)$  and $Y_j(y)=\theta_{\gamma_2}^j(y)$. 
On the probability space $([0,1] \times [0,1], \nu_{\gamma_1} \otimes \nu_{\gamma_2})$ 
the two sequences $(X_i)_{i \geq 1}$ and 
$(Y_j)_{j \geq 1}$ are 
independent. 

\smallskip

Note that, if $\gamma_1=\gamma_2$, the hypothesis
$H_0: \pi=1/2$  is satisfied and the random variables $T_n$ defined in \eqref{Tn} converge in  
distribution to ${\mathcal N}(0,1)$. 

\smallskip

We shall  study the behavior 
of $T_n$  for   different choices of $\gamma_1, \gamma_2, a_n$ and $b_n$. As explained in 
Section \ref{sec2}, the quantities $a_n$ and $b_n$ may be chosen by analyzing the 
graph of the estimated auto-covariances $\hat \gamma_X(k)$ and $\hat \gamma_Y(\ell)$. 

\smallskip

As in Subsection \ref{secar}, we shall  compute $T_n$ for different choices of $(n,m)$ (from $(150,100)$ to 
$(750,500)$ with $n/m=1.5$). We estimate the three quantities 
$\mathrm {Var}(T_n)$, ${\mathbb P}(T_n>1.645)$  and 
 ${\mathbb P}(|T_n|>1.96)$  via a classical Monte-Carlo procedure, by averaging
 over $N=2000$ independent trials.
 
 \smallskip
 
  If $\gamma_1= \gamma_2$,  we say that
${\mathbb P}(T_n>1.645)$ is level 1 (the level of the 
one-sided test)  and ${\mathbb P}(|T_n|>1.96)$ is level 2
(the level of the two-sided test). If 
  $a_n, b_n$ are well chosen, the estimated variance
 should be close to 1 and the estimated levels 1 and 2 should be close to $0.05$. 
 
 \smallskip
 
 If $\gamma_1\neq  \gamma_2$, we say that
${\mathbb P}(T_n>1.645)$ is power 1 (the power of the 
one-sided test)  and ${\mathbb P}(|T_n|>1.96)$ is power
2 (the power of the two-sided test). 

\smallskip

 Note that we 
cannot simulate exactly  $(X_i)_{i \geq 1}$ and 
$(Y_j)_{j \geq 1}$, since  the distribution $\nu_{1/4}$ is unknown.
We start by picking $X_1$ and $Y_1$ independently over $[0, 0.05]$ (in order to start near the neutral fixed point 0) and 
then generate $X_i=\theta_{\gamma_1}^{i-1}(X_1)$ and $Y_j=\theta_{\gamma_2}^{j-1}(Y_1)$ for $i>1, j>1$. 
We shall see that the tests are robust to this (slight) lack of stationarity.

\medskip

$\bullet$ Case $\gamma_1=\gamma_2=1/4$ and  $a_n=0, b_n=0$ (no correction):

\medskip

\begin{center}
\begin{tabular}{|c|c|c|c|c|c|}
\hline
\quad \quad \quad \quad \quad \quad $n ; m$  
&    150 ; 100  & 300 ; 200 & 450 ; 300 & 600 ; 400 & 750 ; 500 \\
\hline 
Estimated variance &5.172  & 4.658  &  4.858 & 4.764 &  4.751 \\
\hline
Estimated level 1 &0.204 &0.196 & 0.19& 0.206 & 0.218\\
\hline
Estimated level 2 &0.358 &0.359 & 0.361 & 0.361& 0.37 \\
\hline
\end{tabular}
\end{center}

\medskip

Here again, since $a_n=b_n=0$, we compute 
$T_n$ as if both series $(X_i)_{i \geq 1}$ and $(Y_j)_{j \geq 1}$ were uncorrelated. 
It leads to a disastrous result, with an estimated variance around 4.8,
an estimated level 1 around $0.2$, and an estimated level 2 around 0.36. This means that the-one sided or two-sided tests build with $T_n$ 
will reject the null hypothesis much too often.

\medskip

\begin{figure}[htb]
\centerline{
\includegraphics[width=11.5 cm, height=6.5 cm]{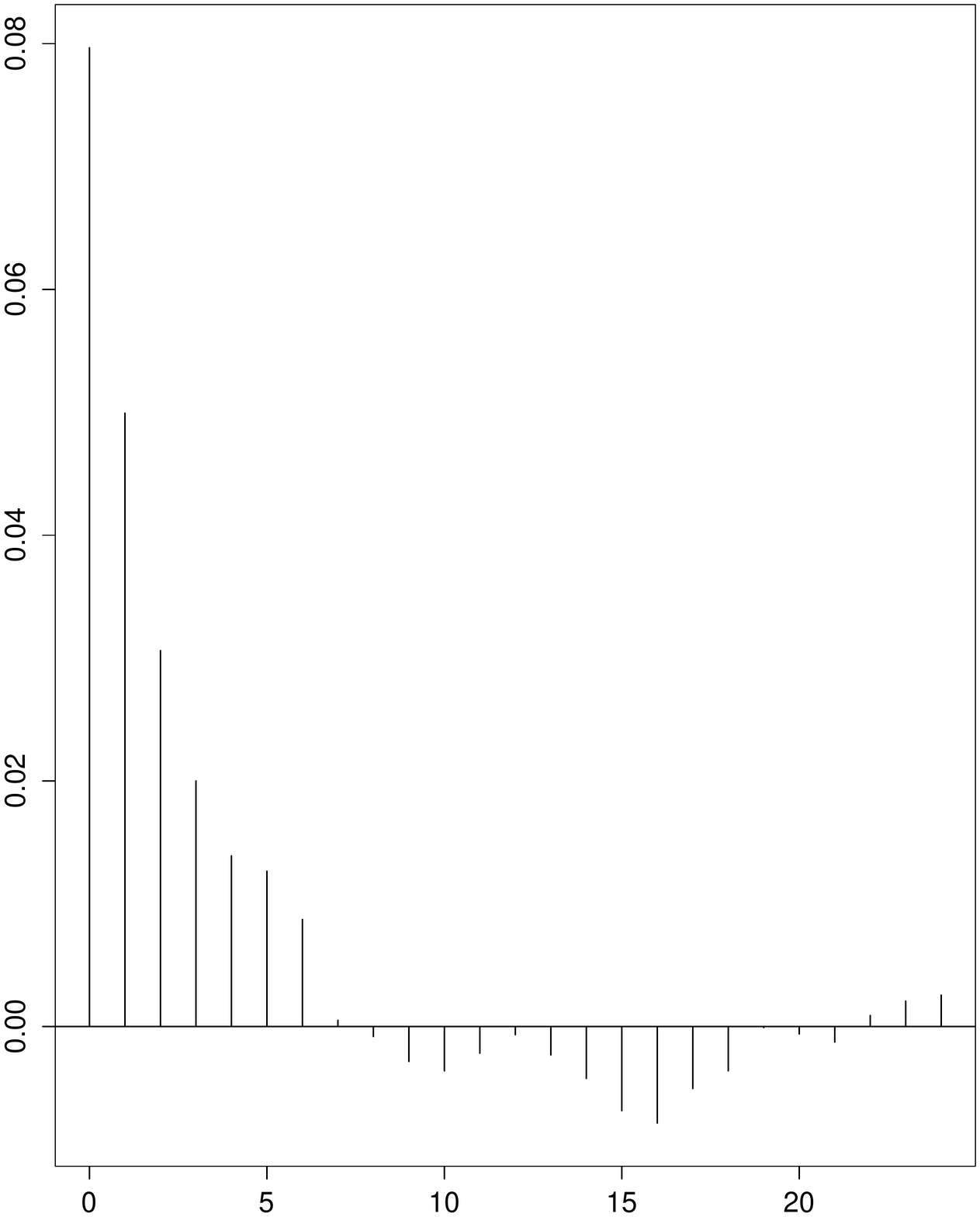}}
\centerline{
\includegraphics[width=11.5 cm, height=6.5 cm]{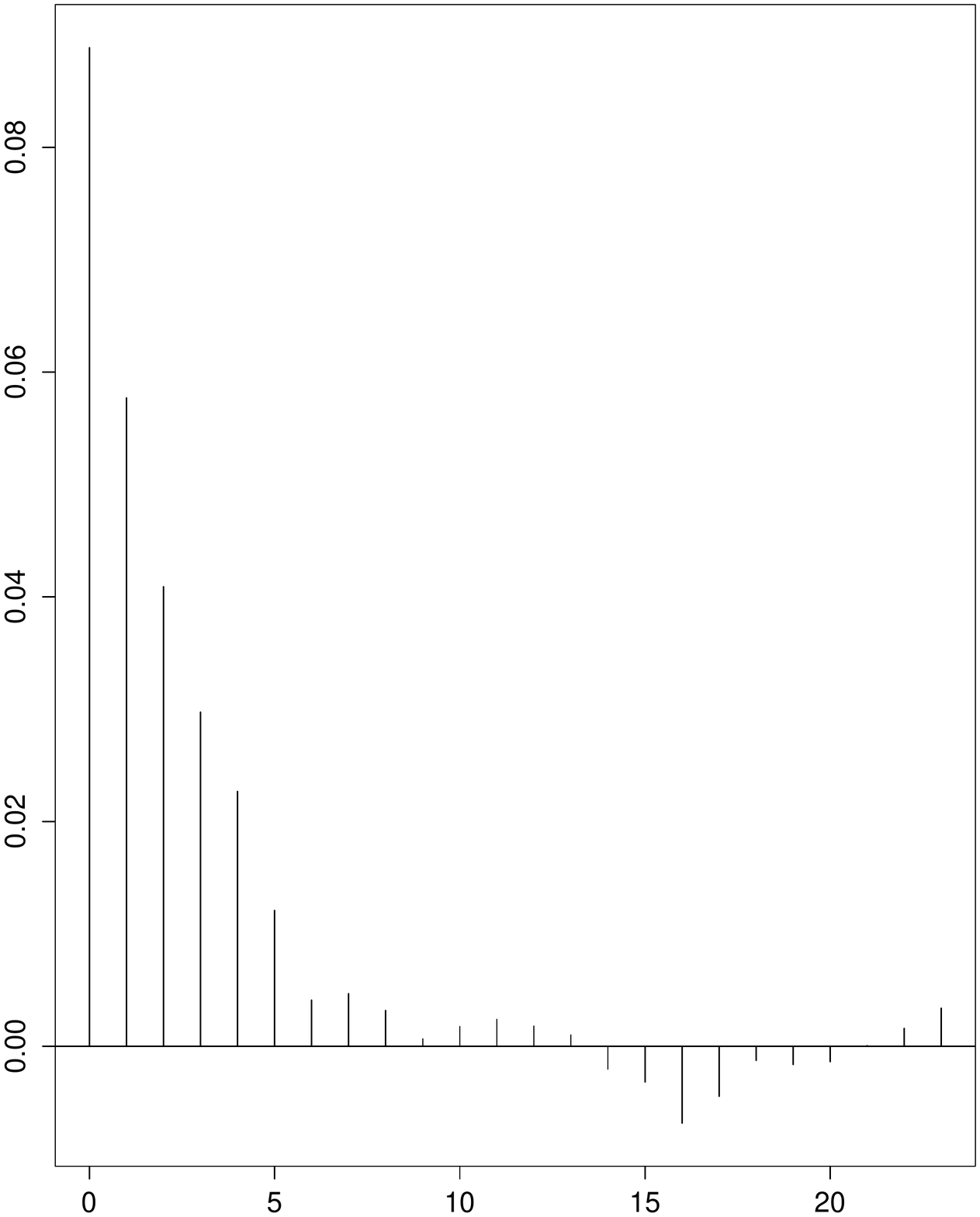}}
\caption{Estimation of the $\gamma_X(k)$'s (top) and 
$\gamma_Y(\ell)$'s (bottom)
for Example 2, 
with   $\gamma_1=\gamma_2=1/4$, $n=300$  and $m=200$.}
\label{fig:covsd}
\end{figure}

$\bullet$ Case  $\gamma_1=\gamma_2=1/4$ and $a_n=4, b_n=4$:

\medskip

\begin{center}
\begin{tabular}{|c|c|c|c|c|c|}
\hline
\quad \quad \quad \quad \quad \quad $n ; m$  
&    150 ; 100  & 300 ; 200 & 450 ; 300 & 600 ; 400 & 750 ; 500 \\
\hline 
Estimated variance &1.352  & 1.308  &  1.329 & 1.324 &  1.28 \\
\hline
Estimated level 1 &0.057 &0.0645 & 0.062 & 0.069 &  0.067\\
\hline
Estimated level 2 &0.089 &0.082 & 0.087 & 0.087& 0.075 \\
\hline
\end{tabular}
\end{center}

\medskip

As suggested by the graphs of the estimated auto-covariances (see Figure 2), the 
choice   $a_n=4, b_n=4$ should give a reasonable estimation of $V$. We see that 
the procedure is much better than if no correction is made: the estimated variance 
is now around $1.3$, the estimated level 1 is around $0.065$ and the estimated level 2 
is around $0.085$.

\medskip

$\bullet$  Case $ \gamma_1=\gamma_2=1/4$ and  $a_n=5, b_n=5$:

\medskip

\begin{center}
\begin{tabular}{|c|c|c|c|c|c|}
\hline
\quad \quad \quad \quad \quad \quad $n ; m$  
&    150 ; 100  & 300 ; 200 & 450 ; 300 & 600 ; 400 & 750 ; 500 \\
\hline 
Estimated variance & 1.349  & 1.249  &  1.164 & 1.181 &  1.185 \\
\hline
Estimated level 1 & 0.057 &0.062 & 0.055 & 0.055 &  0.063\\
\hline
Estimated level 2 & 0.084 &0.077 & 0.068 & 0.071& 0.07  \\
\hline
\end{tabular}
\end{center}

\medskip

For $a_n=5, b_n=5$, the results seem slightly better than for $a_n=4, b_n=4$.
Even for large samples, the estimated level 2 is still larger than the expected level:
around 0.07 instead of 0.05 (but $0.07$ has to be compared to $0.36$ which is the estimated level 2 without correction). Recall that the coefficients $\alpha_{2, \bf X}(n)$ 
and $\alpha_{2, \bf Y}(n)$ are exactly of order $n^{-1/3}$, which is quite slow.
The procedure can certainly be improved for larger $n,m$ by estimating more covariance
terms (in accordance with Proposition \ref{varest}).

\medskip 

$\bullet$ Case $\gamma_1=0.25, \gamma_2=0.1$ and $a_n=5, b_n=4$. 

\medskip

We first estimate the quantity $\pi = {\mathbb P}(X_1< Y_1)$ via 
a realization of 
$$
  \frac {1}{nm} \sum_{i=1}^n \sum_{j=1}^m {\mathbf 1}_{X_i <Y_j} \, .
$$
For $n=m=30000$, this gives an estimation  around $0.529$ for $\pi$, meaning that
$Y$ weakly dominates $X$. Let us see how  the tests based on $T_n$ can see
this  departure from  $ H_0:  \pi=1/2$.

Based on the estimated covariances, a reasonable choice is $a_n=5, b_n=4$.
This choice leads to the following estimated powers (based on $N=2000$ independent
trials).

\medskip

\begin{center}
\begin{tabular}{|c|c|c|c|c|c|}
\hline
\quad \quad \quad \quad \quad \quad $n ; m$  
&    150 ; 100  & 300 ; 200 & 450 ; 300 & 600 ; 400 & 750 ; 500 \\
\hline 
Estimated variance  &  1.24 & 1.2 & 1.2 &  1.23&  1.117\\
\hline
Estimated power 1 &  0.11 & 0.165 & 0.185 &  0.226&  0.241\\
\hline
Estimated power 2 &  0.091 &0.12 & 0.13  &  0.158& 0.17 \\
\hline
\end{tabular}
\end{center}

\medskip

As one can see, the powers of the one sided and two-sided tests are always greater than $0.05$, as 
expected. Still as expected, these powers increase with the size of the samples.
For $n=750, m=500$, the power of the one-sided test is around $0.24$
and that of the two-sided test is around $0.17$. This has to be compared 
with the estimated levels when $T_n$ is centered on 0.529 instead of 0.5 (recall that 
0.529 is the previous 
estimation of $\pi$ on very large samples): in that case, for $n=750, m=500$ and $a_n=5, b_n=4$, 
one gets 
an estimated level 1 = 0.059, and an estimated level 2 = 0.068, which are both close to 
0.05 as expected.

\subsection{Intermittents maps: the case of two adjacent samples}

Here, we shall illustrate the result of Subsection \ref{adj}.
As in Subsection \ref{Sec:IM}, we first simulate $X_1$ according to the uniform distribution over   $[0, 0.05]$,  and 
then generate $X_i=\theta_{\gamma}^{i-1}(Z_1)$, 
but now for the indexes $i \in \{2, \ldots, n+m\}$. 
For $j \in \{n+1, \ldots, n+m \}$, we take 
$Y_j=f(X_j)$ for a given monotonic function $f$.

\smallskip

In our first experiment, we simply take $f=$Id, so that $Y_j=X_j$. Hence 
the assumption $H_0: \pi =1/2$ is satisfied.

\smallskip

For the second experiment, we take $f(x)=x^{4/5}$, so that
$f(X_1)>X_1$ almost surely. Hence $\pi > 1/2$. 

\smallskip

Let $T_n$  be defined in \eqref{Tnbis}.
We follow the same scheme as in Subsection \ref{Sec:IM}, 
and we estimate 
the three quantities 
$\mathrm {Var}(T_n)$, ${\mathbb P}(T_n>1.645)$  and 
 ${\mathbb P}(|T_n|>1.96)$  via a classical Monte-Carlo procedure, by averaging
 over $N=2000$ independent trials.

\medskip

$\bullet$  Case $\gamma=1/4$, $f=$Id and  $a_n=0, b_n=0$ (no correction):

\medskip

\begin{center}
\begin{tabular}{|c|c|c|c|c|c|}
\hline
\quad \quad \quad \quad \quad \quad $n ; m$  
&    150 ; 100  & 300 ; 200 & 450 ; 300 & 600 ; 400 & 750 ; 500 \\
\hline 
Estimated variance & 4.851  & 4.908  &  4.551 & 4.926 &  4.85 \\
\hline
Estimated level 1 & 0.294 &0.276 & 0.252 & 0.256 &  0.251\\
\hline
Estimated level 2 & 0.368 &0.373 & 0.354 & 0.366& 0.362 \\
\hline
\end{tabular}
\end{center}

\medskip

Here, as in the corresponding two-sample case (see the table for $\gamma_1=\gamma_2= 1/4$), the uncorrected 
tests lead to a disastrous result, with an estimated level 1 around
0.25 and an estimated level 2 around 0.36.

\medskip

$\bullet$  Case $\gamma=1/4$, $f=$Id and  $a_n=5, b_n=5$:

\medskip

\begin{center}
\begin{tabular}{|c|c|c|c|c|c|}
\hline
\quad \quad \quad \quad \quad \quad $n ; m$  
&    150 ; 100  & 300 ; 200 & 450 ; 300 & 600 ; 400 & 750 ; 500 \\
\hline 
Estimated variance & 1.369  & 1.302  &  1.272 & 1.205 &  1.171 \\
\hline
Estimated level 1 & 0.114 &0.097 & 0.086 & 0.081 &  0.075\\
\hline
Estimated level 2 & 0.109 &0.096 & 0.079 & 0.072& 0.066  \\
\hline
\end{tabular}
\end{center}

\medskip

The corrected tests with  the choice $a_n=b_n=5$ give a much
better result than if no correction is made. One can observe, 
however, that the estimated levels are still too high ($\geq 8 \%$) for moderately large samples (up to $n=450, m=300$). The results seem less good than in the corresponding two-sample
case (see the table for $\gamma_1=\gamma_2= 1/4$),  which is of course not surprising.

\medskip

$\bullet$  Case $\gamma=1/4$, $f(x)=x^{4/5}$ and  $a_n=5, b_n=5$:

\medskip

\begin{center}
\begin{tabular}{|c|c|c|c|c|c|}
\hline
\quad \quad \quad \quad \quad \quad $n ; m$  
&    150 ; 100  & 300 ; 200 & 450 ; 300 & 600 ; 400 & 750 ; 500 \\
\hline 
Estimated variance & 1.596  & 1.448  &  1.334 & 1.311 &  1.451 \\
\hline
Estimated power 1 & 0.321 &0.385 & 0.451 & 0.49 &  0.548\\
\hline
Estimated power 2 & 0.242 &0.293 & 0.34 & 0.383& 0.442  \\
\hline
\end{tabular}
\end{center}

\medskip

As one can see, the powers of the one sided and two-sided tests are always greater than $0.05$, as 
expected. Still as expected, these powers increase with the size of the samples.
For $n=750, m=500$, the power of the one-sided test is around $0.55$
and that of the two-sided test is around $0.44$.

\bigskip

\noindent  {\bf \large{Acknowledgments.}} We thank Herold Dehling for helpful comments on a preliminary version of the paper.

\end{document}